\documentclass{article}
\usepackage[utf8]{inputenc}
\usepackage{geometry}
\usepackage{graphicx}
\graphicspath{ {./images/} }
\usepackage{amsmath}
\usepackage{amssymb}
\usepackage{listings}
\usepackage{xcolor}
\usepackage{hyperref}
\usepackage{authblk}
\usepackage{setspace}
\title{Tensegrity Structures and Nonlinear Algebra}
\author{Alex Heaton }

\begin{document}

\begin{center}
    \begin{Large} Nonlinear algebra with tensegrity frameworks \end{Large}\\
    Alexander Heaton\\
    March 26, 2020
\end{center}

\begin{abstract}
    In this paper, we discuss tensegrity from the perspective of nonlinear algebra in a manner accessible to undergraduates. We compute explicit examples and include the SAGE and Julia code so that readers can continue their own experiments and computations. The entire framework is a natural extension of linear equations of equilibrium, but to describe the space of solutions will require (nonlinear) polynomials. In our examples, minors of a structured matrix determine the singular locus of the algebraic variety of interest. At these singular points, more interesting phenomena can occur, which we investigate in the context of the tensegrity 3-prism, our running example. Tools from algebraic geometry, commutative algebra, semidefinite programming, and numerical algebraic geometry will be used. Although at first it is all linear algebra, the examples will motivate the study of systems of polynomial equations. In particular, we will see the importance of varieties cut out by determinants of matrices.
\end{abstract}
\begin{singlespace}
    \ref{introduction} Introduction\\
    \ref{section:preliminary-setup} Preliminary Setup\\
    \ref{section:understanding-rigidity-matrix} Intuitive understanding of the rigidity matrix\\
    \ref{section:tangent-space-vector-bundle} Tangent spaces and sections of vector bundles\\
    \ref{section:rigidity-prestress-stability} Rigidity, prestress rigidity, and the singular locus\\
    \ref{section:numerical-algebraic-geometry} Numerical algebraic geometry\\
    \ref{section:adjacent-minors} Gr\"obner bases and primary decomposition\\
\end{singlespace}
\section{Introduction}\label{introduction}
In 1948, the artist and sculptor Kenneth Snelson created a surprisingly stable structure from rigid bars and almost invisible wires \cite{S2012}. It looked like it should collapse, but it didn't. Buckminster Fuller took up the cause and spread similar ideas across the world \cite{F1975}. These structures are fascinating to look at, simply because their existence seems an impossibility. Bars appear to float in midair, and yet cables and bars together can create a remarkably rigid object complete with both \textit{tension} and structural \textit{integrity}.
\begin{figure}[!htb]
\minipage{0.3\textwidth}
  \includegraphics[width=\linewidth]{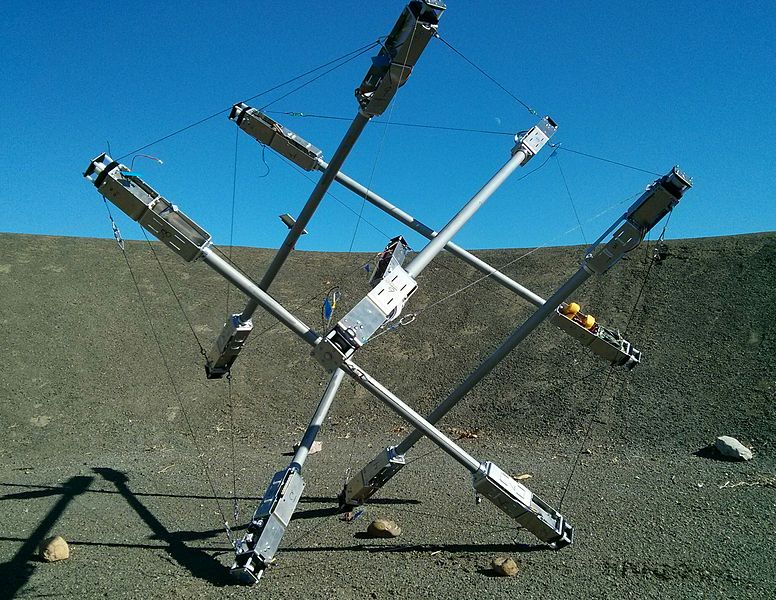}
  \caption{NASA has a robot named Super Ball Bot, inspired by tensegrity structures, adapted for landing on other planets. \cite{superballbotWiki} }\label{figure:NASA-superballbot}
\endminipage\hfill
\minipage{0.26\textwidth}
  \includegraphics[width=\linewidth]{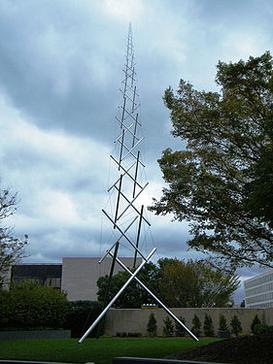}
  \caption{Snelson's Needle Tower, located in Washington D.C. \cite{needletowerWiki}}\label{figure:Snelson-needletower}
\endminipage\hfill
\minipage{0.3\textwidth}%
  \includegraphics[width=\linewidth]{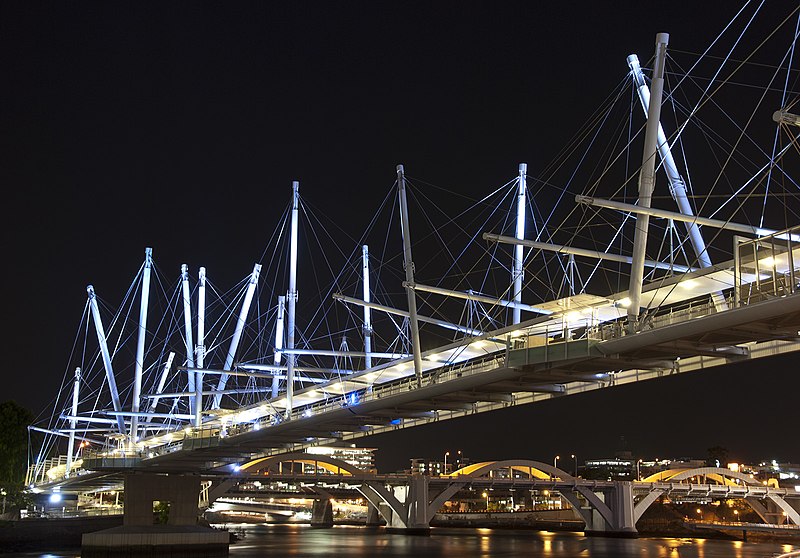}
  \caption{Kurilpa Bridge in Brisbane, Australia, is the world's largest tensegrity bridge. \cite{kurilpabridgeWiki}}\label{figure:Kurilpa-bridge}
\endminipage
\end{figure}

\begin{figure}[!htb]
    \centering
    \includegraphics[width=0.55\textwidth]{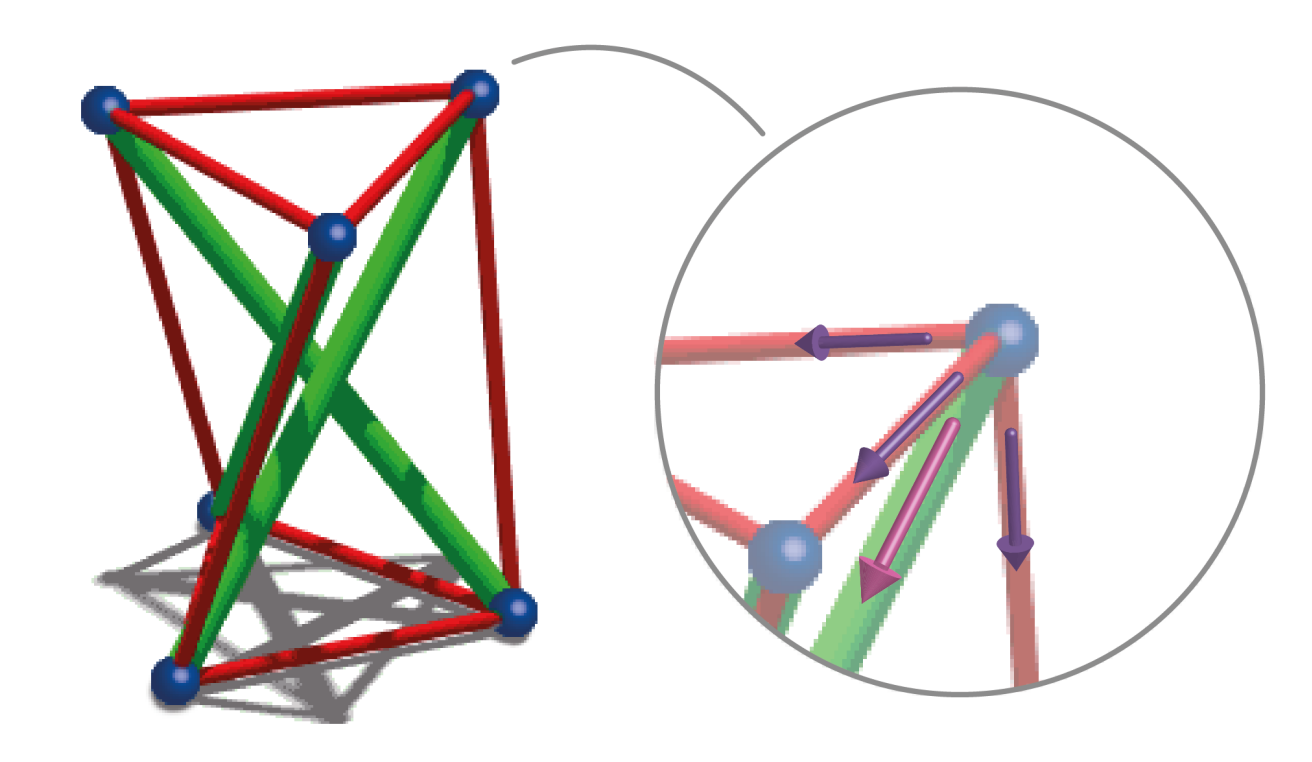}
    \caption{Tensegrity 3-prism, our illustration inspired by \cite{3prismWiki}}
    \label{figure:3-prism-zoom}
\end{figure}

Our running example in this paper will be the \textit{tensegrity 3-prism} pictured in Figure \ref{figure:3-prism-zoom}. This structure consists of 6 nodes, 9 cables, and 3 bars. The cables are stretched tight. Like a child building a toy telephone out of a string and two cups, our cables are useless unless they are under tension. The 3 bars are able to sustain both stretching and compression, remaining the same length.
Think about one of the corners. All the red cables are stretched, in tension, and so exert forces along their individual directions. If you could feel the three blue vectors shown in Figure \ref{figure:3-prism-zoom}, the overall force would be the vector sum, and hence look something like the pink force arrow. But the green bar in our tensegrity 3-prism is \textit{compressed}. It wants to expand, and so exerts a force \textit{exactly and perfectly} opposite to the pink arrow. In fact, the precise locations of the 6 blue nodes of our tensegrity 3-prism are delicate. Even small changes in the coordinates of the nodes could ruin this perfect balancing of forces! If you search for videos of people building tensegrity frameworks, you will get a sense for their tricky behavior.

\section{Preliminary setup}\label{section:preliminary-setup}

We will discuss \textit{bar frameworks} and \textit{tensegrity frameworks} in this paper. Both concepts begin with a combinatorial graph on $n$ nodes $[n]:=\{1,2,\dots,n\}$ and $m$ edges $E \subset \binom{[n]}{2}$. We treat edges as two element subsets of $[n]$, denoting an edge as $\{i,j\} \in E$ or $ij \in E$ for short. For the 3-prism of Figure \ref{figure:3-prism-zoom} we have nodes $[6]=\{1,2,3,4,5,6\}$ and edges $E = \left\{ 12, 13, 14, 15, 23, 25, 26, 34, 36, 45, 46, 56 \right\}$, so that $m = 12$. An \textit{initial configuration} is a map $p:[n] \to \mathbb{R}^d$ usually for $d=2$ or $d=3$. Such a map may be given by listing the coordinates of each node. For the 3-prism embedded as in Figure \ref{figure:3-prism-zoom} we have
\begin{equation}\label{equation:3prism-embedding-p0}
    p = \left[ \begin{array}{ccc}
        p_{11} & p_{12} & p_{13}\\ p_{21} & p_{22} & p_{23}\\ p_{31} & p_{32} & p_{33}\\ p_{41} & p_{42} & p_{43}\\ p_{51} & p_{52} & p_{53}\\ p_{61} & p_{62} & p_{63}
    \end{array} \right] = \left[ \begin{array}{ccc}
        1 & 0 & 0\\ -\frac{1}{2} & \frac{\sqrt{3}}{2} & 0\\  -\frac{1}{2} & -\frac{\sqrt{3}}{2} & 0\\ -\frac{\sqrt{3}}{2} & -\frac{1}{2} & 3\\ \frac{\sqrt{3}}{2} & -\frac{1}{2} & 3\\ 0 & 1 & 3
    \end{array} \right].
\end{equation}
We also abuse notation and think of $p$ as a tuple $(p_{ik}) \in \mathbb{R}^{nd}$ where the $k$th dimensional component of node $i$ is the real number $p_{ik}$. In the computer code accompanying this paper \cite{heaton-tensegrity-trusses-github} we further abuse notation and think of $p$ as an $n \times d$ matrix of real numbers, as above. For an embedding $p:[n] \to \mathbb{R}^d$ the squared edge length $\ell_{ij}^2$ of the edge $ij \in E$ was drawn on a wall for Donald Duck in \textit{Mathemagic land} roughly as
\begin{equation*}
    \ell_{ij}^2 = \sum_{k=1}^d \left(p_{ik} - p_{jk}\right)^2.
\end{equation*}
We partition the edge set (disjointly) into sets of \textit{bars, cables} and \textit{struts} as in $E = B \cup C \cup S$. For every bar $ij \in B$ we introduce an equation $g_{ij}(x) = 0$ given by the polynomial
\begin{align*}
    g_{ij}(x) &= \sum_{k=1}^d ( x_{ik} - x_{jk} )^2 - \ell_{ij}^2\\
     &= \sum_{k=1}^d ( x_{ik} - x_{jk} )^2 - (p_{ik} - p_{jk})^2.
\end{align*}
Cables are allowed to shorten, so we replace $g_{ij}(x)=0$ by $g_{ij}(x) \leq 0$ for $ij \in C$. Struts are allowed to lengthen, so we have $g_{ij}(x) \geq 0$ for $ij \in S$. These equations and inequalities allow us to consider all possible configurations $(x_{ik}) \in \mathbb{R}^{nd}$ for a given graph $([n],E)$ which satisfy the \textit{member constraints}
\begin{equation}\label{equation:member-constraints}
    \begin{array}{ccc}
        g_{ij}(x)=0 & & ij \in B  \\
        g_{ij}(x)\leq 0 & & ij \in C  \\
        g_{ij}(x)\geq 0 & & ij \in S.  \\
    \end{array}
\end{equation}
If $E=B$ so that all edges are bars, we say we have a \textit{bar framework}, and the \textit{valid configurations} $(x_{ik}) \in \mathbb{R}^{nd}$ satisfying (\ref{equation:member-constraints}) form an algebraic variety, since they are solutions to a finite list of polynomial equations. If $E \neq B$ then there are cables or struts and we have a \textit{tensegrity framework}. In that case, the valid configurations form a semi-algebraic set, since they are solutions to finitely many polynomial equations and inequalities. Usually we are given an embedding $p:[n] \to \mathbb{R}^d$, so that we immediately start with one solution $(x_{ik}) = (p_{ik}) \in \mathbb{R}^{nd}$ to the member constraints (\ref{equation:member-constraints}). We write $g(p)=0$. But you could also choose the real numbers $\ell_{ij}^2$ in some other way, and then you might start without any solutions $(x_{ik}) \in \mathbb{R}^{nd}$ to (\ref{equation:member-constraints}). The remainder of this paper explores some available tools to find and understand these solutions. In particular, we will explore tools from algebraic geometry and commutative algebra like Gr\"obner bases and primary decomposition, but also tools from semidefinite programming and numerical algebraic geometry. Every matrix, equation, and formula, along with Figures $5,6,8,11,12,13,14,15,16,17,20,21,24$ can be reproduced from \texttt{SAGE} or \texttt{Julia} code freely available at \cite{heaton-tensegrity-trusses-github}.

\begin{figure}[!htb]
    \centering
    \includegraphics[width=0.25\textwidth]{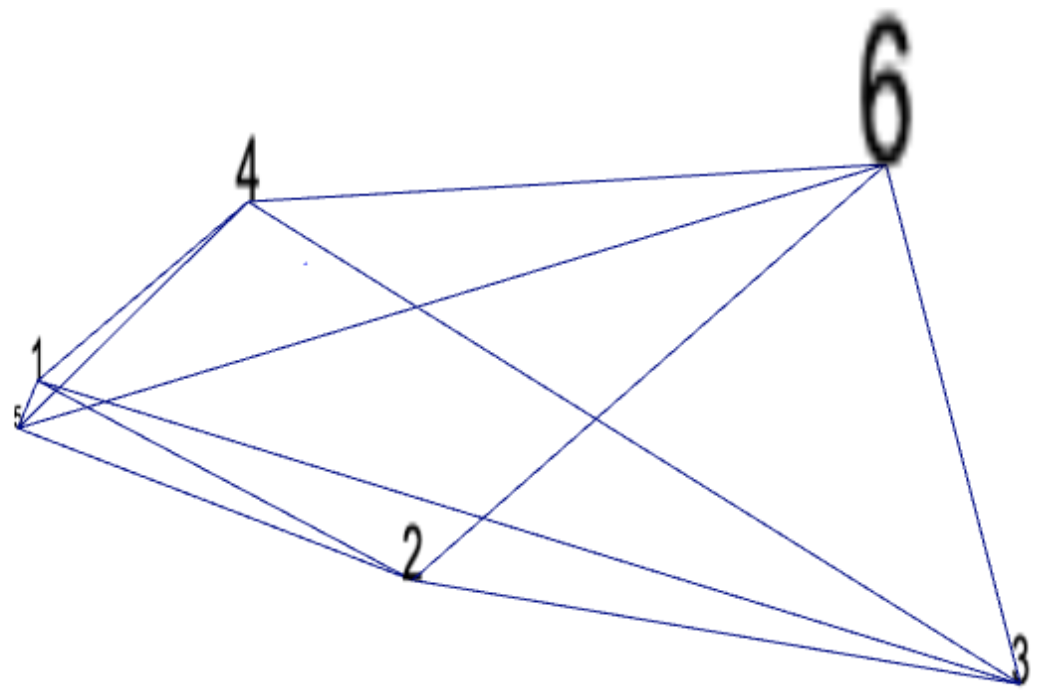}
    \includegraphics[width=0.13\textwidth]{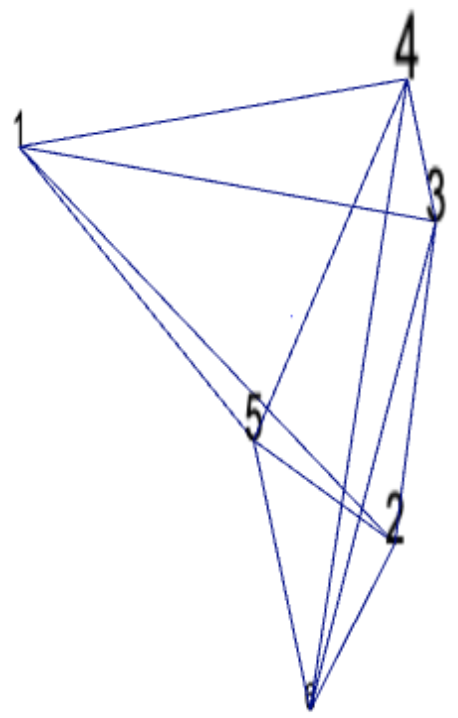}
    \includegraphics[width=0.25\textwidth]{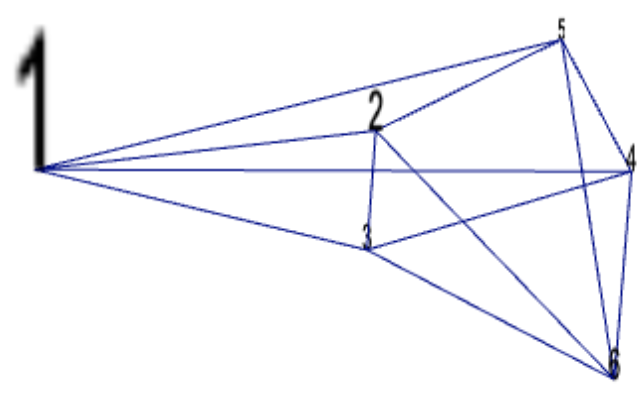}
    \includegraphics[width=0.20\textwidth]{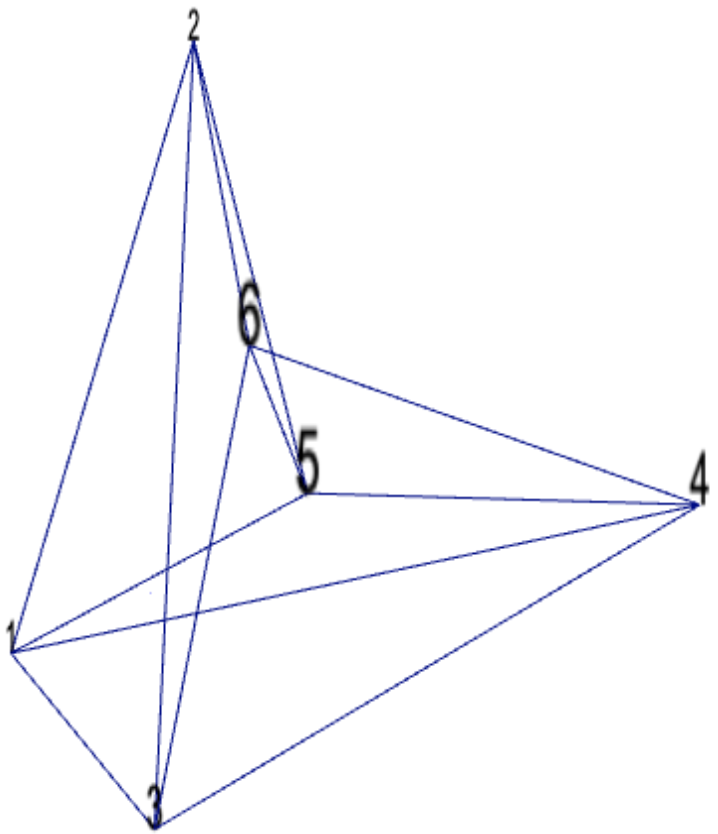}
    \caption{Some random 3-prisms in $\mathbb{R}^3$}
    \label{figure:random}
\end{figure}

To begin, we describe some geometry in the case of bar frameworks $E=B$, since there we have only equations. Since there are $m$ edges we have $m$ polynomials in $nd$ variables. In fact, the geometry is easier if we allow our variables to take on complex number values. Don't let this put you off. All our familiar real-valued solutions, including $(p_{ik}) \in \mathbb{R}^{nd}$, are still solutions to the complex system. It is helpful to write this system of equations briefly as $g(x) = 0$ where $x=(x_{ik})$ and $0 \in \mathbb{C}^m$. Here $g$ is a polynomial map $g:\mathbb{C}^{nd} \to \mathbb{C}^m$. Writing it this way reminds us there are $nd$ variables and $m$ equations. We also write $g$ as a column vector $g = [ g_1(x), g_2(x), \dots, g_m(x) ]^T$ whose entries are polynomials. Here we label the equations by $1,2,\dots,m$ instead of $ij \in E$. For the 3-prism with embedding $p$ as in (\ref{equation:3prism-embedding-p0}) we have 
\begin{equation}\label{equation:member-constraints-3prism} \footnotesize
    g(x) = \left[ \begin{array}{c}
    {\left(x_{11} - x_{21}\right)}^{2} + {\left(x_{12} - x_{22}\right)}^{2} + {\left(x_{13} - x_{23}\right)}^{2} - 3\\
{\left(x_{11} - x_{31}\right)}^{2} + {\left(x_{12} - x_{32}\right)}^{2} + {\left(x_{13} - x_{33}\right)}^{2} - 3\\
{\left(x_{11} - x_{41}\right)}^{2} + {\left(x_{12} - x_{42}\right)}^{2} + {\left(x_{13} - x_{43}\right)}^{2} - \frac{1}{4} \, {\left(\sqrt{3} + 2\right)}^{2} - \frac{37}{4}\\
{\left(x_{11} - x_{51}\right)}^{2} + {\left(x_{12} - x_{52}\right)}^{2} + {\left(x_{13} - x_{53}\right)}^{2} - \frac{1}{4} \, {\left(\sqrt{3} - 2\right)}^{2} - \frac{37}{4}\\
{\left(x_{21} - x_{31}\right)}^{2} + {\left(x_{22} - x_{32}\right)}^{2} + {\left(x_{23} - x_{33}\right)}^{2} - 3\\
{\left(x_{21} - x_{51}\right)}^{2} + {\left(x_{22} - x_{52}\right)}^{2} + {\left(x_{23} - x_{53}\right)}^{2} - \frac{1}{2} \, {\left(\sqrt{3} + 1\right)}^{2} - 9\\
{\left(x_{21} - x_{61}\right)}^{2} + {\left(x_{22} - x_{62}\right)}^{2} + {\left(x_{23} - x_{63}\right)}^{2} - \frac{1}{4} \, {\left(\sqrt{3} - 2\right)}^{2} - \frac{37}{4}\\
{\left(x_{31} - x_{41}\right)}^{2} + {\left(x_{32} - x_{42}\right)}^{2} + {\left(x_{33} - x_{43}\right)}^{2} - \frac{1}{2} \, {\left(\sqrt{3} - 1\right)}^{2} - 9\\
{\left(x_{31} - x_{61}\right)}^{2} + {\left(x_{32} - x_{62}\right)}^{2} + {\left(x_{33} - x_{63}\right)}^{2} - \frac{1}{4} \, {\left(\sqrt{3} + 2\right)}^{2} - \frac{37}{4}\\
{\left(x_{41} - x_{51}\right)}^{2} + {\left(x_{42} - x_{52}\right)}^{2} + {\left(x_{43} - x_{53}\right)}^{2} - 3\\
{\left(x_{41} - x_{61}\right)}^{2} + {\left(x_{42} - x_{62}\right)}^{2} + {\left(x_{43} - x_{63}\right)}^{2} - 3\\
{\left(x_{51} - x_{61}\right)}^{2} + {\left(x_{52} - x_{62}\right)}^{2} + {\left(x_{53} - x_{63}\right)}^{2} - 3
    \end{array} \right] = \left[ \begin{array}{c}
        0\\0\\0\\0\\0\\0\\0\\0\\0\\0\\0\\0
    \end{array} \right]. \normalsize
\end{equation}
The algebraic variety $V(g)$ and the real algebraic variety $V_{\mathbb{R}}(g)$ are the sets
\begin{align*}
    V(g) &:= \left\{ x \in \mathbb{C}^{nd} \, : \, g(x) = 0 \right\}\\
    V_{\mathbb{R}}(g) &:= V(g) \cap \mathbb{R}^{nd}.
\end{align*}
These sets are complicated, but what can we say about them? The two most basic invariants are \textit{dimension} and \textit{degree}. Precise algebraic or geometric definitions can be found in many places including \cite{CoxLittleOshea2015IdealsVarietiesAlgorithmsTEXT, SommeseWampler2005numericalSolutionofSystemsofPolynomialsTEXT}. For us, the first step in nonlinear algebra is to recall linear algebra. To that end we examine the Jacobian matrix $dg$ whose entries are polynomials. In the $i$th row and $j$th column we find the partial derivative of the $i$th equation with respect to the $j$th variable. Since our equations and variables are already labelled in many ways, this can get notationally confusing. An example is best. For the tensegrity 3-prism we have
\begin{equation}\label{equation:3-prism-dg-matrix} \footnotesize
    \frac{1}{2}dg = \begin{array}{c}
        \left(\begin{array}{rrrrrrr}
x_{11} - x_{21} & x_{12} - x_{22} & x_{13} - x_{23} & -x_{11} + x_{21} & -x_{12} + x_{22} & -x_{13} + x_{23} & 0 \\
x_{11} - x_{31} & x_{12} - x_{32} & x_{13} - x_{33} & 0 & 0 & 0 & -x_{11} + x_{31} \\
x_{11} - x_{41} & x_{12} - x_{42} & x_{13} - x_{43} & 0 & 0 & 0 & 0 \\
x_{11} - x_{51} & x_{12} - x_{52} & x_{13} - x_{53} & 0 & 0 & 0 & 0 \\
0 & 0 & 0 & x_{21} - x_{31} & x_{22} - x_{32} & x_{23} - x_{33} & -x_{21} + x_{31} \\
0 & 0 & 0 & x_{21} - x_{51} & x_{22} - x_{52} & x_{23} - x_{53} & 0 \\
0 & 0 & 0 & x_{21} - x_{61} & x_{22} - x_{62} & x_{23} - x_{63} & 0 \\
0 & 0 & 0 & 0 & 0 & 0 & x_{31} - x_{41} \\
0 & 0 & 0 & 0 & 0 & 0 & x_{31} - x_{61} \\
0 & 0 & 0 & 0 & 0 & 0 & 0 \\
0 & 0 & 0 & 0 & 0 & 0 & 0 \\
0 & 0 & 0 & 0 & 0 & 0 & 0
\end{array}\right. \cdots \\
  \\
        \cdots \begin{array}{rrrrrr}
0 & 0 & 0 & 0 & 0 & 0 \\
-x_{12} + x_{32} & -x_{13} + x_{33} & 0 & 0 & 0 & 0 \\
0 & 0 & -x_{11} + x_{41} & -x_{12} + x_{42} & -x_{13} + x_{43} & 0 \\
0 & 0 & 0 & 0 & 0 & -x_{11} + x_{51} \\
-x_{22} + x_{32} & -x_{23} + x_{33} & 0 & 0 & 0 & 0 \\
0 & 0 & 0 & 0 & 0 & -x_{21} + x_{51} \\
0 & 0 & 0 & 0 & 0 & 0 \\
x_{32} - x_{42} & x_{33} - x_{43} & -x_{31} + x_{41} & -x_{32} + x_{42} & -x_{33} + x_{43} & 0 \\
x_{32} - x_{62} & x_{33} - x_{63} & 0 & 0 & 0 & 0 \\
0 & 0 & x_{41} - x_{51} & x_{42} - x_{52} & x_{43} - x_{53} & -x_{41} + x_{51} \\
0 & 0 & x_{41} - x_{61} & x_{42} - x_{62} & x_{43} - x_{63} & 0 \\
0 & 0 & 0 & 0 & 0 & x_{51} - x_{61}
\end{array} \cdots \\
          \\
        \cdots \left. \begin{array}{rrrrr}
0 & 0 & 0 & 0 & 0 \\
0 & 0 & 0 & 0 & 0 \\
0 & 0 & 0 & 0 & 0 \\
-x_{12} + x_{52} & -x_{13} + x_{53} & 0 & 0 & 0 \\
0 & 0 & 0 & 0 & 0 \\
-x_{22} + x_{52} & -x_{23} + x_{53} & 0 & 0 & 0 \\
0 & 0 & -x_{21} + x_{61} & -x_{22} + x_{62} & -x_{23} + x_{63} \\
0 & 0 & 0 & 0 & 0 \\
0 & 0 & -x_{31} + x_{61} & -x_{32} + x_{62} & -x_{33} + x_{63} \\
-x_{42} + x_{52} & -x_{43} + x_{53} & 0 & 0 & 0 \\
0 & 0 & -x_{41} + x_{61} & -x_{42} + x_{62} & -x_{43} + x_{63} \\
x_{52} - x_{62} & x_{53} - x_{63} & -x_{51} + x_{61} & -x_{52} + x_{62} & -x_{53} + x_{63}
\end{array}\right)
    \end{array} \normalsize
\end{equation}
We apologize if this $12 \times 18$ matrix hurt your eyes, but it is that important. In the computer code provided at \cite{heaton-tensegrity-trusses-github} you can generate random values for each variable $x_{ik}$ and calculate the rank of $dg|_q$ evaluated at this random point $q \in \mathbb{R}^{nd}$, as well as plot the resulting configurations as in Figure \ref{figure:random}.
You will consistently find $dg|_q$ has rank 12. Since there are 18 columns, this leaves a 6-dimensional nullspace. It turns out these 6 dimensions correspond to the Euclidean group of rigid motions, having 3 translational and 3 rotational degrees of freedom. Since there is no 7th dimension of nullspace, the configurations of Figure \ref{figure:random} are \textit{infinitesimally rigid} and hence \textit{rigid}. We explain this further in Section \ref{section:rigidity-prestress-stability}. For now, we will hint at the role the \textit{generic rank} of $dg$ plays in deciding the \textit{local dimension} of the set $V(g)$ at our given solution $p$. As is further explained in \cite{SommeseWampler2005numericalSolutionofSystemsofPolynomialsTEXT, WamplerHauensteinSommese2011mechanismMobility} we have the following inequalities:
\begin{equation}\label{equation:bounding-local-dimension}
    \binom{d+1}{2} \leq \text{corank }(dg|_q) \leq \text{dim}_{p} V(g) \leq \text{corank }(dg|_p).
\end{equation}
Here, the corank is the dimension of the nullspace of the matrix $dg|_q$ evaluated at a \textit{generic} point $q \in \mathbb{C}^{nd}$ or of $dg|_p$ evaluated at our given point $p = (p_{ik})$. This can be calculated by Gaussian elimination, for example. We postpone an explanation until Section \ref{section:tangent-space-vector-bundle}, hoping that at least we have motivated the study of this matrix $dg$. In the next Section \ref{section:understanding-rigidity-matrix} we hope to give an intuitive explanation of the role $dg$ plays in the context of bar and tensegrity frameworks.

\section{Intuitive understanding of the rigidity matrix}\label{section:understanding-rigidity-matrix}
    
\begin{figure}[!htb]
    \centering
    \includegraphics[width=0.9\textwidth]{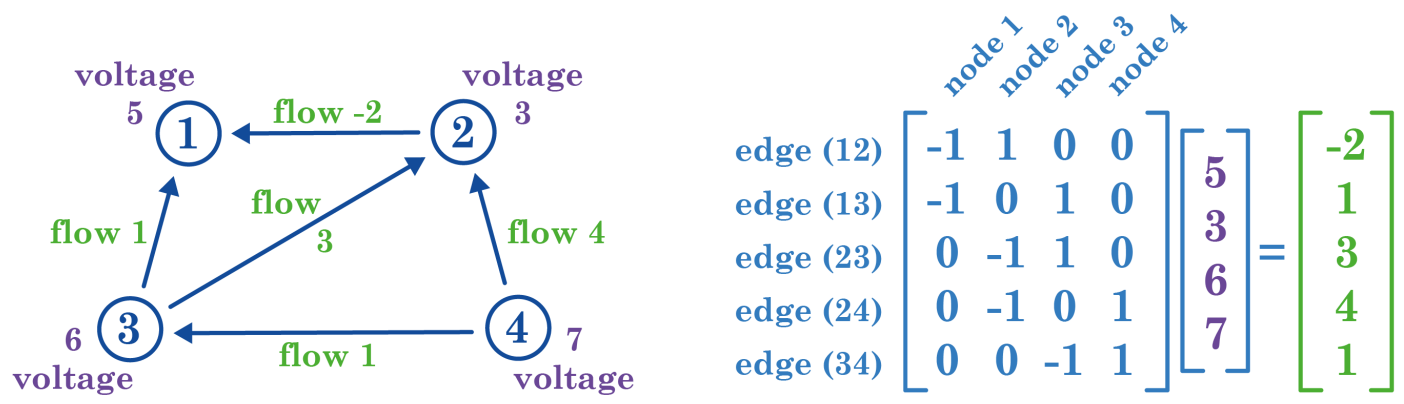}
    \caption{An electrical network and $Av$}
    \label{figure:electrical-network}
\end{figure}

Towards understanding the matrix $dg$, we will make a brief analogy with electrical networks, following the Monthly article \cite{Strang1989PatternsLinearAlgebraAMM} or also \cite{Strang1988frameworkForEquilibriumEquationsIncidenceMatrix}. We imagine voltages at each node, driving currents of electicity to flow through the edges. Current flows from high voltage to low voltage, so what we really need to know are the \textit{voltage differences}. Figure \ref{figure:electrical-network} includes a matrix gathering the voltage differences for us.
Following \cite{Strang1988frameworkForEquilibriumEquationsIncidenceMatrix, Strang1989PatternsLinearAlgebraAMM}, we call this matrix the \textit{incidence matrix} $A$ and we will create an analogous matrix for tensegrity frameworks. Each row comes from an edge, and each row has a $1$ and a $-1$. At this point we note that the voltage differences in the vector $Av$ could be scaled by a diagonal matrix of positive \textit{conductances} (not pictured above) giving the current that would flow along each edge of our graph. These currents can be recorded in the vector $C A v$. What if we want to record the net flow at each individual node? A miracle of nature occurs. Because each column of our incidence matrix corresponds to a specific node in our graph, we can add up all the currents flowing in and out of each node simply by multiplying with $A^T$. In this case we obtain:
\begin{equation*}
    \begin{bmatrix} -1&-1&0&0&0\\1&0&-1&-1&0\\0&1&1&0&-1\\0&0&0&1&1 \end{bmatrix} \begin{bmatrix}-2\\1\\3\\4\\1 \end{bmatrix} = \begin{bmatrix}1\\-9\\3\\5 \end{bmatrix}
\end{equation*}
The first entry tells us that $1$ unit of current would \textit{leave} node 1, while the second entry tells us $9$ units of current would \textit{enter} node 2. Analogous interpretations exist for the other entries.

Of course, to start analyzing actual electrical networks, we would need to use a few laws (Kirchhoff and Ohm). Since our goal is to briefly mention electrical networks as an entry-way into understanding $dg$, we will not venture farther in this direction. For an introduction, we recommend \cite{S1986}. The main take-away is that we created a matrix $A$ which took voltages at nodes, turning their voltage differences into \textit{flows} along the edges. $A$ also played another role as $A^T$. $A^T$ sums the net flows into each particular node. While $A$ went from nodes to edges, $A^T$ went from edges to nodes. The basic equations can be written as follows, where $f$ stands for \textit{external forces or flows}.
\begin{equation*}
    A^T C A v = f
\end{equation*}
In fact, for bar and tensegrity frameworks the same equations will be relevant. There $A$ is called the \textit{rigidity matrix} \cite{ConnellyWhiteley1996secondOrderRigidityAndPrestressStabilityTensegrityFrameworks} while $A^TCA$ is the \textit{stiffness matrix} \cite{Strang2007computationalScienceEngineeringTEXT}. In yet other contexts, this matrix is called the \textit{weighted graph Laplacian}. $A^TCA$ is a positive semi-definite matrix, which can be made positive definite by \textit{grounding a node} in the case of electrical networks, removing the all ones vector from the nullspace (see Equation \ref{equation:eigenpairs}). Similarly, for tensegrity frameworks, we will need to \textit{ground} the $\binom{d+1}{2}$-dimensional group of rigid motions, whose null vectors include the all ones vector spread out across $d$ translations, but also $\binom{d}{2}$ rotations. We will get there soon. First, we build the matrix $A$ by realizing how it should take forces on the nodes into stresses along the edges. Finally, we will understand $A^T$ as summing up the effects of stressed edges on each node. In fact the delicate balance of stresses referenced in the introduction and Figure \ref{figure:3-prism-zoom} will translate to a vector equation $w^T A = 0$ or $A^T w = 0$. Such $w$ are called \textit{self stresses} in \cite{ConnellyWhiteley1996secondOrderRigidityAndPrestressStabilityTensegrityFrameworks}. In electrical networks $w^T A = 0$ represents flow around loops. Euler's topological formula \textit{nodes} $-$ \textit{edges} $+$ \textit{loops} $=1$ comes from subspace dimensions $\text{dim }\mathbb{R}^n - \text{dim }\mathbb{R}^m + \text{dim Null}A^T = \text{dim Null}A$. But let's leave electrical networks and step up to bar frameworks with $d \geq 2$.
\begin{center}
    How to generalize $-1$ and $1$?
\end{center}
\begin{figure}[!htb]
\minipage{0.25\textwidth}
  \includegraphics[width=\linewidth]{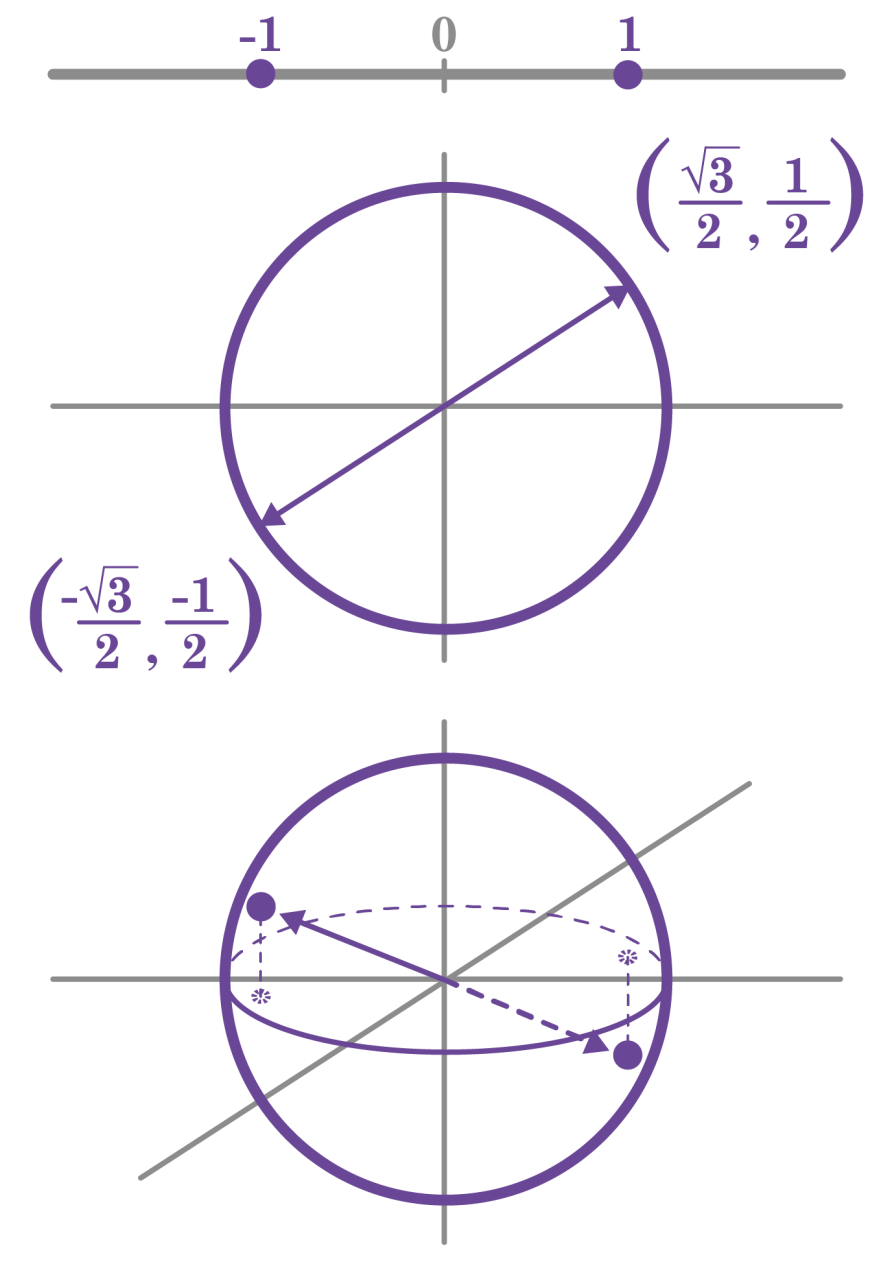}
    \caption{Unit vectors}
    \label{figure:unit-vectors}
\endminipage\hfill
\minipage{0.72\textwidth}
  \includegraphics[width=\linewidth]{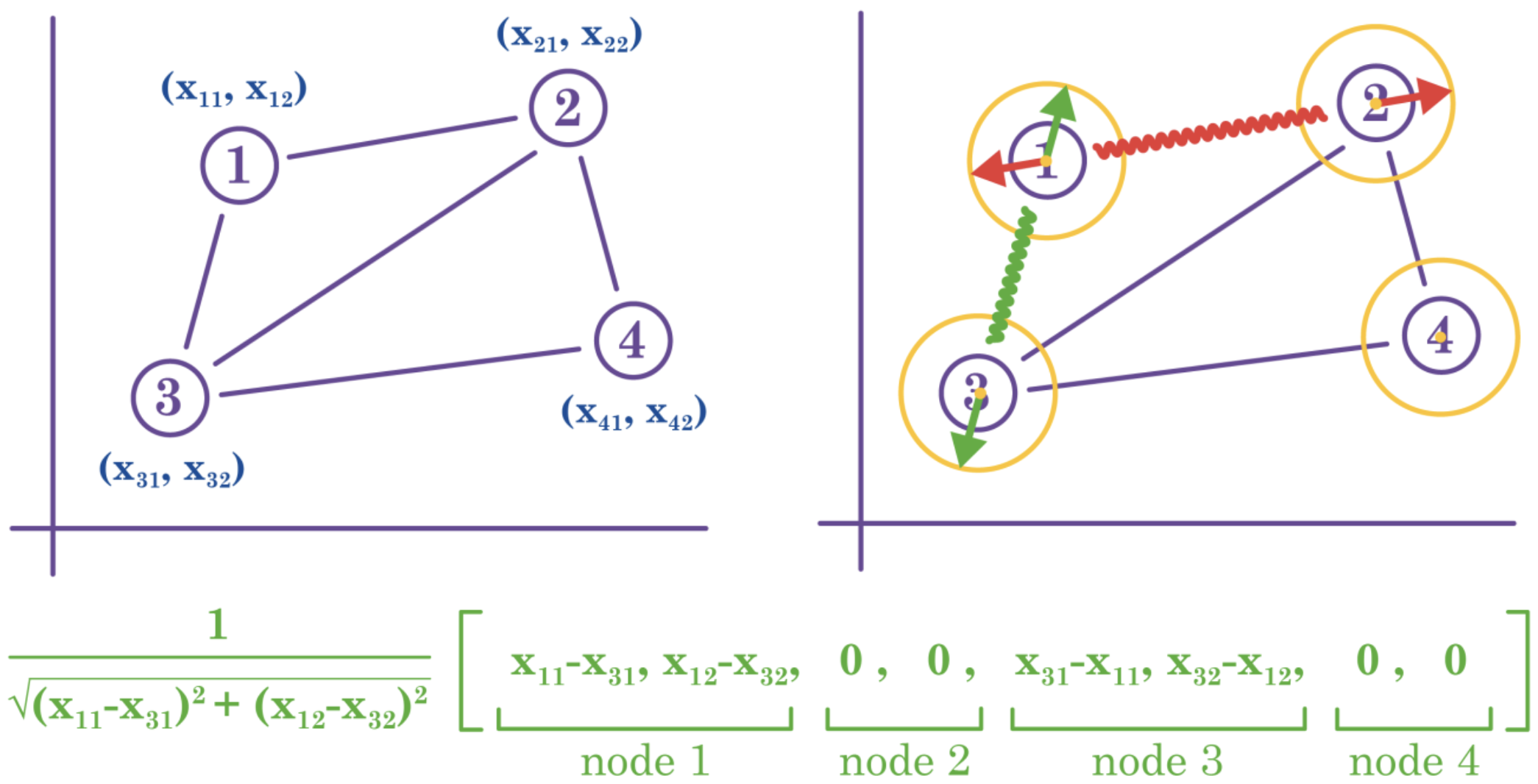}
  \caption{Building the matrix one row per edge}\label{figure:building-rows}
\endminipage
\end{figure}
To build the incidence matrix for an electrical network, it's easiest to build it one row at a time. To build the rigidity matrix $A$ of a tensegrity framework in the plane $\mathbb{R}^2$, or in space $\mathbb{R}^3$, or even in any higher dimension $\mathbb{R}^d$, we will proceed similarly. We build $A$ one row at a time, and almost all entries of our row will be zeros. The only change is to replace $1$ and $-1$ with antipodal points on the relevant sphere of correct dimension (Figure \ref{figure:unit-vectors}). If our framework lives in the plane we take antipodal points on the unit circle. If our framework lives in space, we take points on the colloquial sphere. In $\mathbb{R}^4$ we use unit vectors on the sphere $S^3 \subset \mathbb{R}^4$.

Consider building $A$ for the framework of Figure \ref{figure:building-rows}. To build the row of our matrix corresponding to the red edge between nodes 1 and 2 we would use the red unit vectors in the picture. Although they are drawn on different circles, all the circles are unit circles and the two red vectors are antipodal points on the unit sphere $S^1 \subset \mathbb{R}^2$.
But since a point on the unit circle has two coordinates, now there are two entries of our row corresponding to each node. In $\mathbb{R}^3$ we will have antipodal unit vectors on $S^2 \subset \mathbb{R}^3$, and so each node will need 3 entries in the row of our matrix.

\begin{figure}[!htb]
\minipage{0.44\textwidth}
  \includegraphics[width=\linewidth]{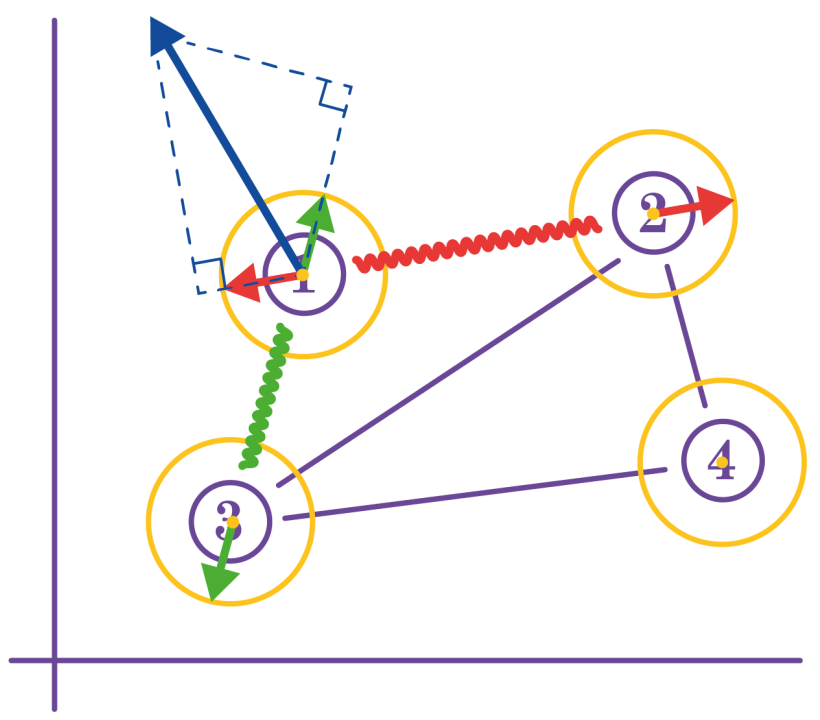}
    \caption{Unsuspecting node 1}
    \label{figure:unsuspecting-force}
\endminipage\hfill
\minipage{0.4\textwidth}
  \includegraphics[width=\linewidth]{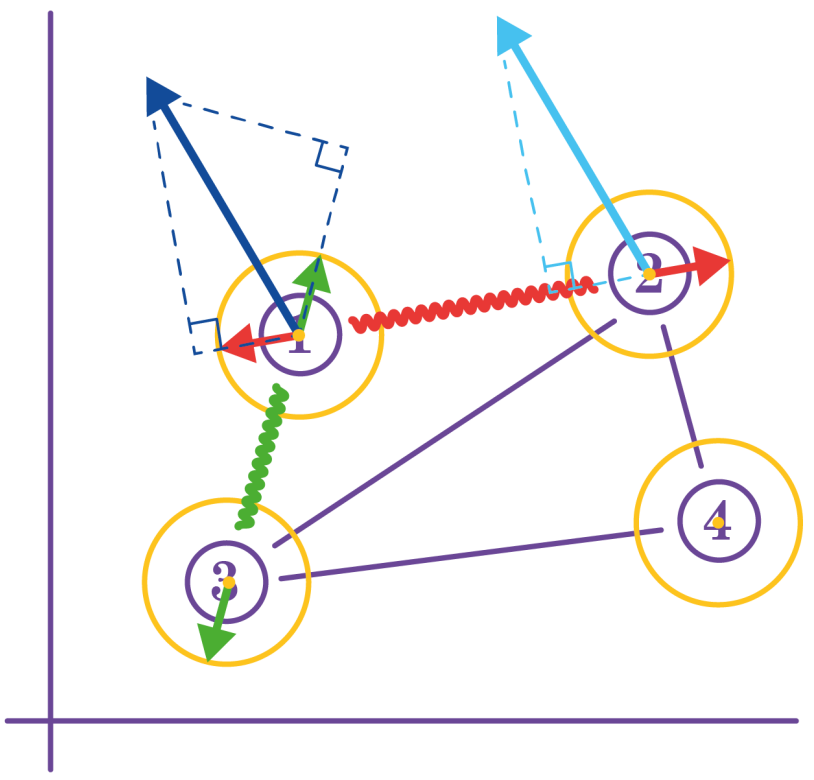}
  \caption{Forces compete on red edge}\label{figure:forces-compete}
\endminipage
\end{figure}

Figure \ref{figure:building-rows} shows the row corresponding to the green edge explicitly. The green edge connects nodes 1 and 3, and so our unit vectors will go in the first two entries, skip 2 entries of zeros, and then fill the next two entries. In this way, you can easily build the matrix $A$ one row/edge at a time. But why is this the case? Why should we generalize the incidence matrix $A$ in this way? Consider the blue force vector of Figure \ref{figure:unsuspecting-force}.
Imagine our poor, unsuspecting framework is standing still. Then, a large blue force vector pulls on node 1. At least in the beginning, this blue force vector will induce tensions in the red and green edges connected to node 1. These tensions are sometimes called internal forces. Treating the edges as springs, Hooke's law tells us the internal force induced on the green edge will be proportional to the dot product of the blue vector with the green unit vector, since we treat this force as giving a displacement, but only count the displacement along the direction of the edge. The internal force induced on the red edge will likewise come from the dot product of the blue force with the red unit vector (Figure \ref{figure:unsuspecting-force}).

Recall for electrical networks, it was the voltage differences that drove current through the edge. Will this happen again? Consider yet another, light-blue force vector, this time on node 2, displayed in Figure \ref{figure:forces-compete}.
Instead of a voltage difference, we have a difference of dot products. Just like $1$ and $-1$ produced a difference in voltages, now our unit vectors (which point in opposite directions) produce a difference in dot products. The light-blue vector has a negative dot product with its red unit vector, causing the red edge to want to compress. The dark-blue vector has a positive dot product with its red unit vector, causing the red edge to want to stretch. Now we see the importance of taking antipodal points on the unit sphere, in order to keep everything straight. Suppose that the dark-blue vector has a larger dot product, then our red edge will \textit{overall} feel a stretching internal force. It's all working out nicely. In fact, it should be clear that the story will be the same in any dimension. Even if we have structures in $\mathbb{R}^5$ feeling $5$-dimensional force vectors, antipodal points on the sphere $S^4 \subset \mathbb{R}^5$ should cover it.

\section{Tangent spaces and sections of vector bundles}\label{section:tangent-space-vector-bundle}

\begin{figure}[!htb]
    \centering
    \includegraphics[width=0.24\textwidth]{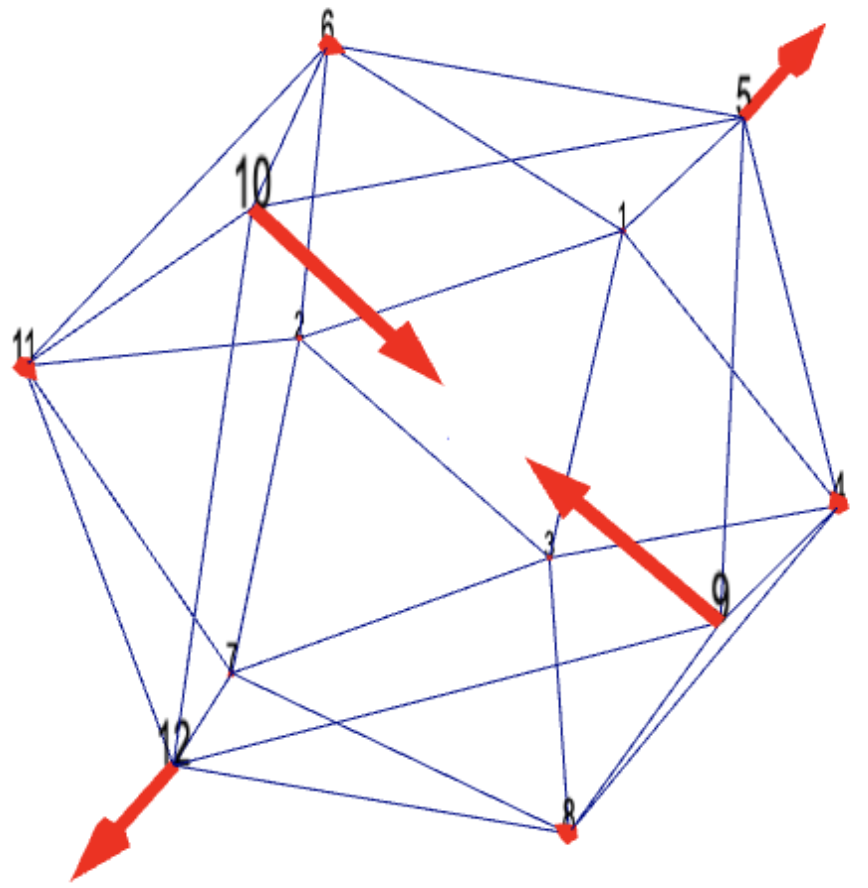}
    \includegraphics[width=0.25\textwidth]{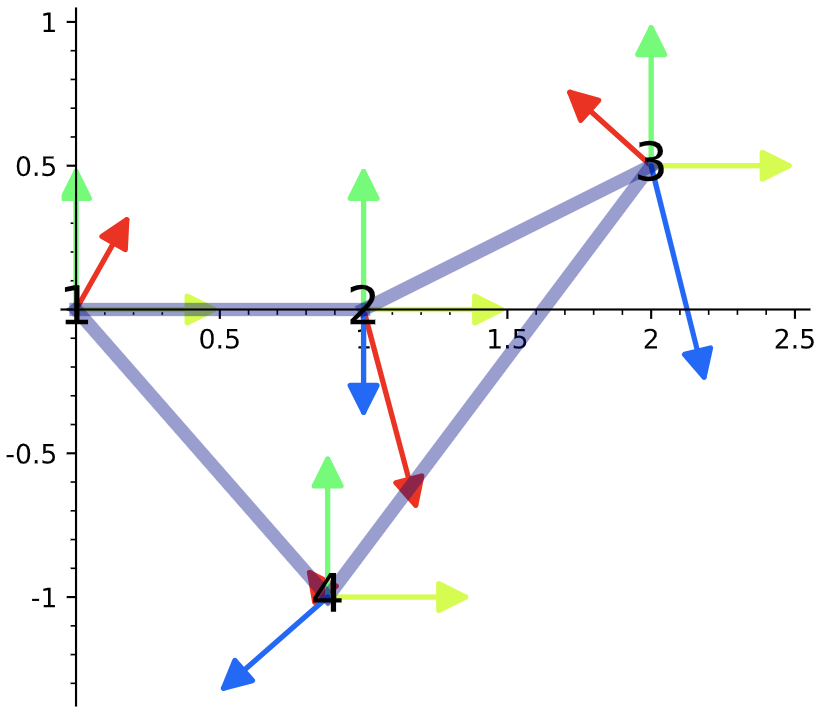}
    \includegraphics[width=0.25\textwidth]{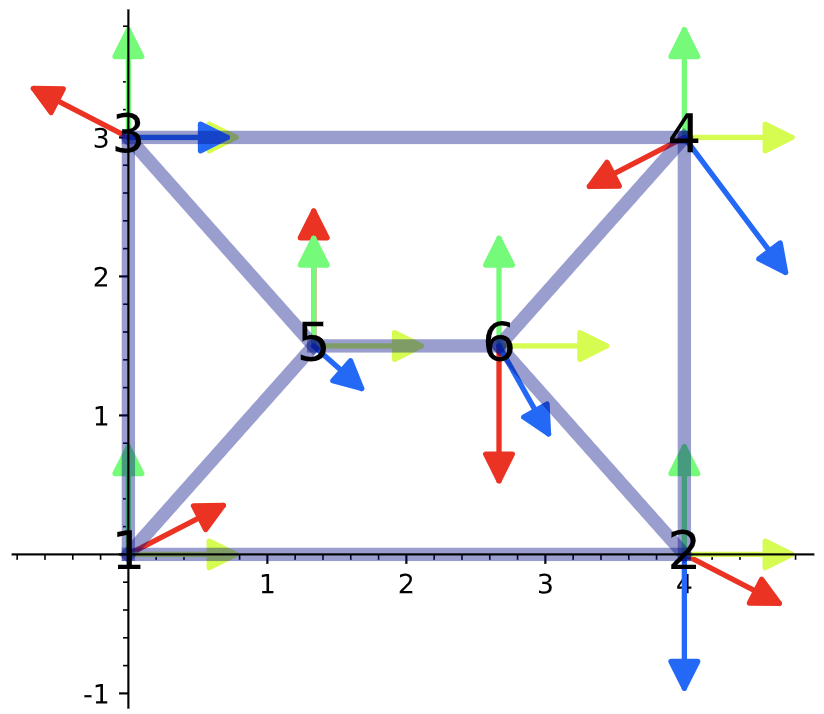}
    \includegraphics[width=0.18\textwidth]{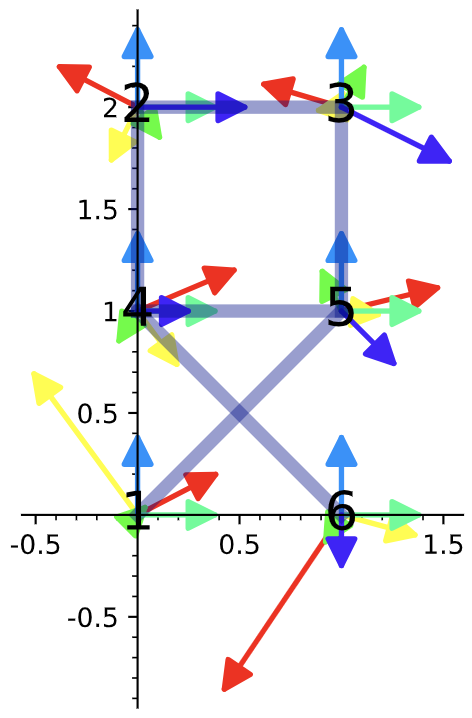}
    \caption{Some infinitesimal mechanisms}
    \label{figure:mechanism-examples}
\end{figure}

You may have noticed the similarity between the Jacobian matrix $dg$ given by partial derivatives and the rigidity matrix $A$ discussed in the previous Section \ref{section:understanding-rigidity-matrix}. Let's take a moment to precisely connect them. To model the edges as springs, we needed anitpodal unit vectors in each row, so we had to normalize, dividing by the square root of $\ell_{ij}^2$, otherwise known as $\ell_{ij}$. This division by a square root occurs in Figure \ref{figure:building-rows}. Because the edges have possibly different lengths, these division factors are different for every row. We collect all the edge lengths $\ell_{ij}$ as entries in a diagonal matrix $L$. Then the relation is
\begin{equation}\label{equation:relating-rigidity-and-jacobian-dg}
    L A = \frac{1}{2} dg.
\end{equation}
$L$ clears the normalizing denominators from each row of $A$ separately, while $\frac{1}{2}$ gets rid of the extra factor $2$ coming from $\frac{d}{dx}(x^2) = 2x$. As hinted in the inequality (\ref{equation:bounding-local-dimension}), understanding the nullspace of $dg$ is helpful. Since $L$ is a full-rank, square, diagonal matrix with positive entries, Equation (\ref{equation:relating-rigidity-and-jacobian-dg}) tells us the nullspace of $A$ is exactly equal to the nullspace of $dg$! This means that we can understand a null vector $dg|_p \cdot v = 0$ for $v \in \mathbb{R}^{nd}$ as an assignment of $n$ smaller $\mathbb{R}^d$ vectors attached to each node in our framework. Thinking of the edges as springs, a null vector $dg|_p \cdot v = 0$ assigns displacements to every node in a way that \textit{does not stretch any spring/edge}. Such vectors are called \textit{infinitesimal mechanisms} in \cite{Strang2007computationalScienceEngineeringTEXT}. The code accompanying this paper automatically plots each null vector in a distinct color, spread across each of the nodes of the framework. A few examples are shown in Figure \ref{figure:mechanism-examples}, but you can easily create your own.
\begin{figure}[!htb]
\minipage{0.25\textwidth}
  \includegraphics[width=\linewidth]{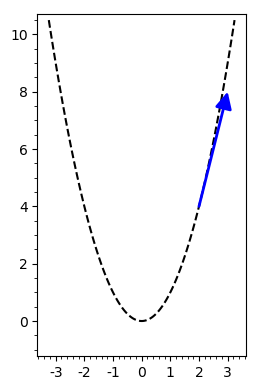}
    \caption{Tangent space}
    \label{figure:tangent-space-parabola}
\endminipage\hfill
\minipage{0.35\textwidth}
  \includegraphics[width=\linewidth]{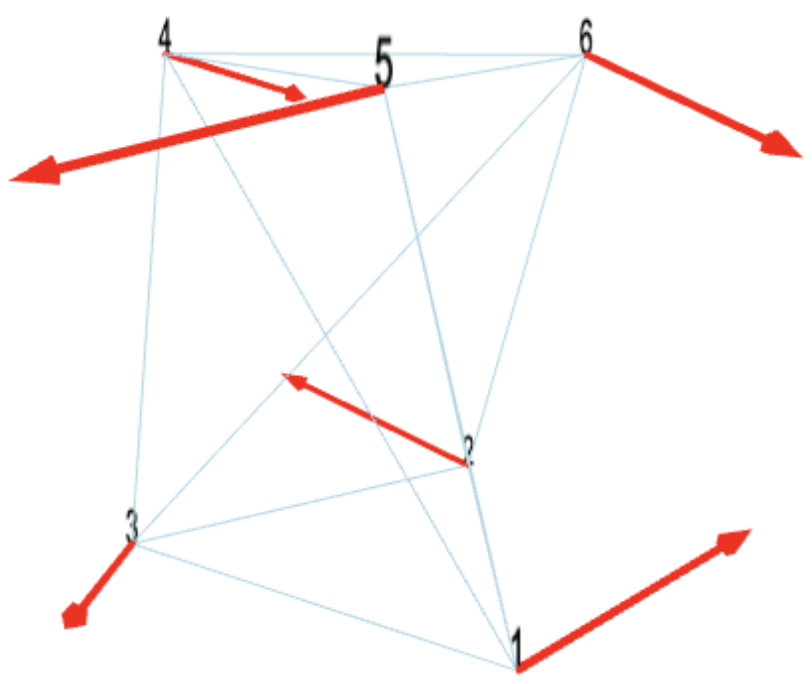}
  \caption{Infinitesimal flexes}\label{figure:3prism-flexes}
\endminipage\hfill
\minipage{0.35\textwidth}
  \includegraphics[width=\linewidth]{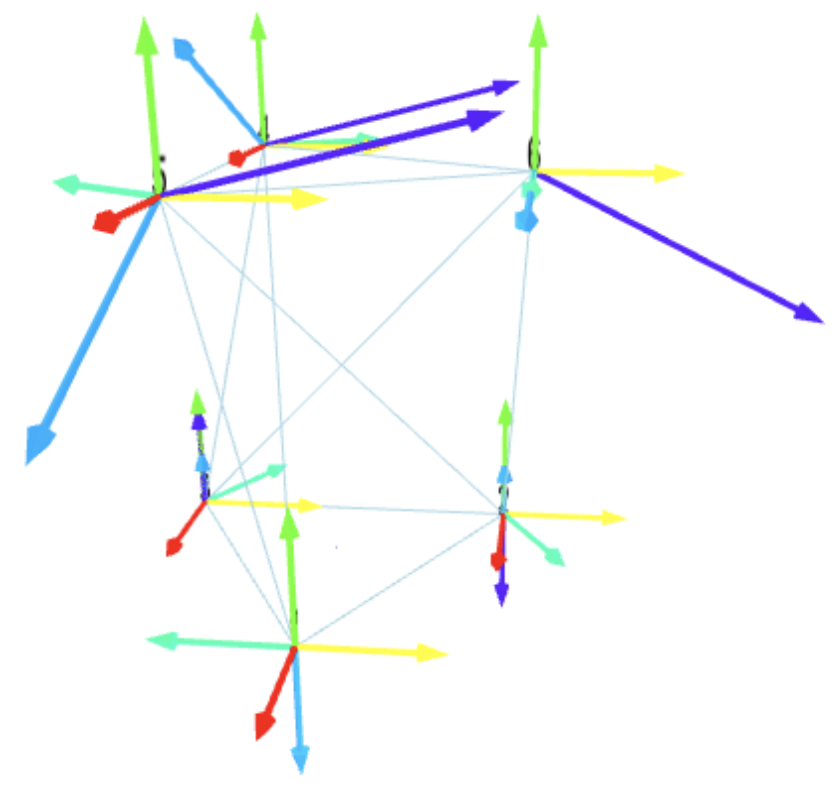}
  \caption{Infinitesimal rigid motions}\label{figure:3prism-rigid-motions}
\endminipage
\end{figure}

Figures \ref{figure:3prism-flexes} and \ref{figure:3prism-rigid-motions} concern the 3-prism, and plot each null vector in a different color. For example, in Figure \ref{figure:3prism-rigid-motions} all the dark blue vectors come from a single null vector $dg|_p \cdot v = 0$ arising from a rotation, while the red, yellow, and green vectors come from the 3 translational rigid motions. The teal and light-blue vectors are the remaining two rotations. Recall from Section \ref{section:preliminary-setup} that a randomly chosen embedding $q:[6] \to \mathbb{R}^3$ for the 3-prism resulted in $\binom{d+1}{2} = 6$ independent null vectors. However, for the (very special) embedding $p$ of Equation (\ref{equation:3prism-embedding-p0}) there is one extra, independent, $7$th null vector! This null vector is plotted in Figure \ref{figure:3prism-flexes}. We emphasize once again that the actual null vectors are column vectors $v \in \mathbb{R}^{18}$ which when multiplying our $12 \times 18$ matrix $dg|_p$ (whose polynomial entries are evaluated at $p$) give the zero vector. For every framework, we have the uniquely determined subspace $R$ of infinitesimal rigid motions, so that we write
\begin{equation}\label{equation:nullspace-decomposition-R-plus-F}
    \text{Null }dg|_p = R \oplus F
\end{equation}
Although a choice of the complementary space $F$ of infinitesimal flexes is not unique, the space $R$ is, and this will be enough for us in Section \ref{section:rigidity-prestress-stability} where we investigate \textit{prestress rigidity}. For general algebraic varieties $V(f)$ we have no way to draw tangent vectors, but for bar frameworks, we can.

About algebraic varieties: Given any polynomial map $f:\mathbb{C}^{N} \to \mathbb{C}^m$ recall that $V(f)$ is the set of all points $x \in \mathbb{C}^N$ such that $f(x) = 0 \in \mathbb{C}^m$, meaning all $m$ polynomials vanish simultaneously. The set $V(f)$ has a geometric structure, including \textit{tangent spaces}, which we describe now. A quick example. If $f:\mathbb{C}^2 \to \mathbb{C}^1$ is given by $f(x) = \left[ f_1(x_1,x_2) = x_2 - x_1^2 \right]$, then $V_{\mathbb{R}}(f) = V(f) \cap \mathbb{R}^2$ is exactly the parabola of Figure \ref{figure:tangent-space-parabola}. The Jacobian is an $m \times N$ matrix $df = \left[ \begin{array}{cc} -2x_1 & 1 \end{array} \right]$ with a one dimensional nullspace at all points. For example, at the point $p = (2,4)$ we can calculate a basis for the nullspace
\begin{equation*}
    \text{Null } df|_p = \text{ span } \left\{ \left[ \begin{array}{c}
        1 \\
        4 
    \end{array} \right] \right\}.
\end{equation*}
In this case, since we have $N=2$ variables, we can plot this vector attached to the point $p=(2,4)$ of the set $V_\mathbb{R}(f)$. For higher dimensional algebraic varieties $N \geq 4$, we cannot visualize $V(f)$ or its tangent spaces very easily. In fact, things get even more interesting.

\begin{figure}[!htb]
\minipage{0.43\textwidth}
  \begin{equation*}
    f(x) = \left[ \begin{array}{c}
        f_1 = -x^{3} + x y\\
        f_2 = -x^{4} + x z\\
        f_3 = x^{7} - x^{5} y - x^{4} z + x^{2} y z
    \end{array} \right]
\end{equation*}
    \caption{A polynomial map}
    \label{figure:a-polynomial-map}
\endminipage\hfill
\minipage{0.48\textwidth}
    \centering
    \includegraphics[width=0.7\linewidth]{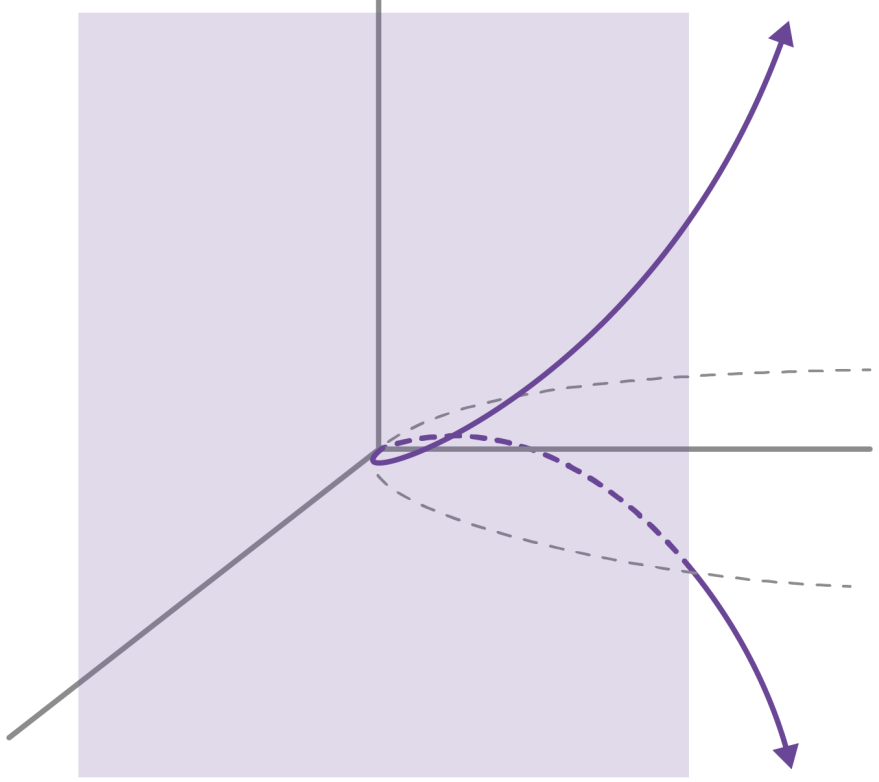}
    \caption{Two irreducible components}\label{figure:twisted-cubic}
\endminipage
\end{figure}

The algebraic variety $V(f)$ for some polynomial map $f:\mathbb{C}^N \to \mathbb{C}^m$ will in general have several \textit{irreducible components}. Consider $f:\mathbb{C}^3 \to \mathbb{C}^3$ given in Figure \ref{figure:a-polynomial-map}, where the $x$ in $f(x)$ actually means $(x,y,z)$.
Then $V(f)$ happens to have 2 irreducible components, call them $X_1$ and $X_2$. $X_1$ is the $x=0$ plane and has dimension $2$, while $X_2$ is the \textit{twisted cubic} and has dimension $1$. Intersected with the reals $V(g) \cap \mathbb{R}^2$ we can see these components in Figure  \ref{figure:twisted-cubic}. We can also look at the Jacobian $df$ and evaluate it at several different points.
\begin{equation*} \footnotesize
    \begin{array}{cccc}
    df & df|_{(0,0,0)} & df|_{(0,5,3)} & df|_{(1,1,1)} \\
     & & & \\
    \left[\begin{array}{rrr}
-3 x^{2} + y & x & 0 \\
-4 x^{3} + z & 0 & x \\
7 x^{6} - 5 x^{4} y - 4 x^{3} z + 2 x y z & -x^{5} + x^{2} z & -x^{4} + x^{2} y
\end{array}\right] & \left[\begin{array}{rrr}
0 & 0 & 0 \\
0 & 0 & 0 \\
0 & 0 & 0
\end{array}\right] & \left[\begin{array}{rrr}
5 & 0 & 0 \\
3 & 0 & 0 \\
0 & 0 & 0
\end{array}\right] & \left[\begin{array}{rrr}
-2 & 1 & 0 \\
-3 & 0 & 1 \\
0 & 0 & 0
\end{array}\right]
    \end{array}
\end{equation*}
Without precise definitions, we see that at \textit{smooth points} $(0,5,3) \in X_1$ and $(1,1,1) \in X_2$ the dimension of the nullspace agrees with the dimension of the component. But there is also the point $(0,0,0) \in V(f)$ whose nullspace is $3$-dimensional. This is a \textit{singular point}. We will discuss this soon.

\begin{figure}[!htb]
    \centering
    \includegraphics[width=0.4\textwidth]{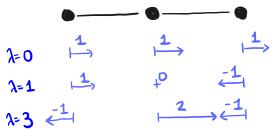}
    \caption{Three independent vibrations}
    \label{figure:molecule-vibrates}
\end{figure}

We have seen that null vectors $df v = 0$ cannot be drawn if the number of variables is more than 3, but they can be drawn for bar frameworks. The reason is that we interpret $v$ as a section of a vector bundle on a finite set. Therefore we cannot resist a remark on \cite[Chapter 3]{Sternberg1994grouptheoryandphysics}, the incidence matrix $A$, and our null vectors. Consider a fictitious molecule with three atoms and two bonds as in Figure \ref{figure:molecule-vibrates}.
Using the incidence matrix $A$ we can create the stiffness matrix $A^T C A$ and find its three \textit{eigenpairs}: multiply each vector to check.
\begin{equation}\label{equation:eigenpairs}
    \begin{bmatrix} 1&0\\-1 & 1\\0&-1 \end{bmatrix} \begin{bmatrix} 1&0\\0&1\end{bmatrix} \begin{bmatrix}1&-1&0\\0&1&-1\end{bmatrix} = \begin{bmatrix}1&-1&0\\-1&2&-1\\0&-1&1\end{bmatrix} \hspace{1cm} \, \, (0,\begin{bmatrix}1\\1\\1\end{bmatrix}) \, \, (1,\begin{bmatrix}1\\0\\-1\end{bmatrix}) \, \, (3,\begin{bmatrix}-1\\2\\-1\end{bmatrix})
\end{equation}
The strongest molecular vibration corresponds to the largest eigenvalue with resonant frequency $\lambda = 3$, which gives an interesting wobbly motion which is also depicted in Figure \ref{figure:molecule-vibrates}. Incidentally, this same eigenbasis is useful for the group $S_3$ which permutes the three basis vectors of $\mathbb{R}^3$. Rewriting the permutation matrices in this eigenbasis, we find that all six permutation matrices simultaneously block-diagonalize. Here is one:
\begin{align*}
    \begin{bmatrix} &1 & \\ & &1\\1& & \end{bmatrix} \mapsto &  \begin{bmatrix} \frac{1}{\sqrt{3}} & 0 & 0\\ 0 & \frac{1}{\sqrt{2}} & 0\\ 0 & 0 & \frac{1}{\sqrt{6}} \end{bmatrix} \begin{bmatrix} 1 & 1 & 1\\ 1 & 0 & -1\\ -1 & 2 & -1 \end{bmatrix} \begin{bmatrix} &1 & \\ & &1\\1& & \end{bmatrix} \begin{bmatrix} 1 & 1 & -1\\1&0&2\\1&-1&-1\end{bmatrix} \begin{bmatrix} \frac{1}{\sqrt{3}} & 0 & 0\\ 0 & \frac{1}{\sqrt{2}} & 0\\ 0 & 0 & \frac{1}{\sqrt{6}} \end{bmatrix}\\
     & = \begin{bmatrix} 3^{-\frac{1}{2}} & 2^{-\frac{1}{2}} & -6^{-\frac{1}{2}}\\ 3^{-\frac{1}{2}} & 0 & 2\cdot 6^{-\frac{1}{2}}\\ 3^{-\frac{1}{2}} & -2^{-\frac{1}{2}} & -6^{-\frac{1}{2}} \end{bmatrix}^T \begin{bmatrix} &1 & \\ & &1\\1& & \end{bmatrix} \begin{bmatrix} 3^{-\frac{1}{2}} & 2^{-\frac{1}{2}} & -6^{-\frac{1}{2}}\\ 3^{-\frac{1}{2}} & 0 & 2\cdot 6^{-\frac{1}{2}}\\ 3^{-\frac{1}{2}} & -2^{-\frac{1}{2}} & -6^{-\frac{1}{2}} \end{bmatrix}\\
     & = \begin{bmatrix}1& & \\ &-\frac{1}{2}& \frac{\sqrt{3}}{2}\\  & -\frac{\sqrt{3}}{2} & -\frac{1}{2} \end{bmatrix}\\
\end{align*}
The process of block-diagonalizing matrices is a fundamental task in the subject of \textit{representation theory}.

This paper is too short to contain a full discussion of the relevance of symmetry and representation theory to these topics. But since the paper's purpose is to inspire future reading and thinking, we point you in one such direction. Very briefly, many representations begin with a symmetry group $G$ acting on a set $X$. This makes each group element into an invertible linear operator on the vector space of all complex-valued functions on the set. $S_3$ can permute the three atoms of our molecule, giving rise to the $3$-dimensional representation we block-diagonalized. Chapter 3 of \cite{Sternberg1994grouptheoryandphysics} discusses how if the tensegrity framework is invariant under the action of a group $G$, then we obtain a representation of $G$ on the space of displacements, sections of a vector bundle on a finite set. Informally, you can think of these as vector-valued functions on the nodes, as in the Figures \ref{figure:mechanism-examples}, \ref{figure:3prism-flexes}, \ref{figure:3prism-rigid-motions}. Since the `molecule' in Figure \ref{figure:molecule-vibrates} has a graph structure, its symmetry group is $S_2$, which permutes the two outer atoms. The stiffness matrix is a self-adjoint operator that commutes with this representation and each eigenvector spans an irreducible representation of $S_2$ (twice the trivial representation, and once the sign representation). In the articles \cite{CT1995,BC1998} the authors study our 3-prism and many other tensegrity frameworks by using symmetry. The dihedral group and other finite subgroups of $SO_3$ can be used to algorithmically produce tensegrity frameworks. In fact, new frameworks were discovered using this method.

\section{Rigidity, prestress rigidity, and the singular locus}\label{section:rigidity-prestress-stability}

\begin{figure}[!htb]
\minipage{0.47\textwidth}
  \includegraphics[width=\linewidth]{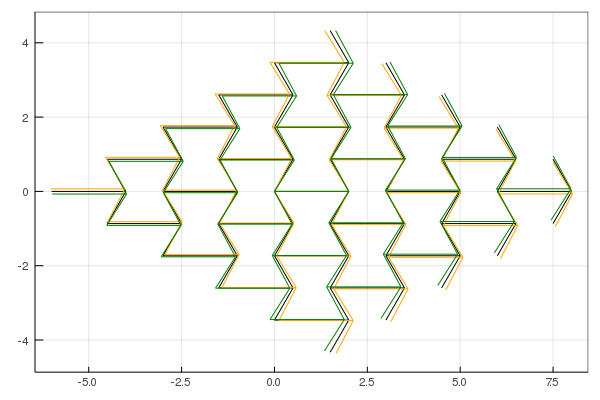}
    \caption{Auxetic honeycomb flexing}
    \label{figure:flexible-honeycomb}
\endminipage\hfill
\minipage{0.43\textwidth}
  \includegraphics[width=\linewidth]{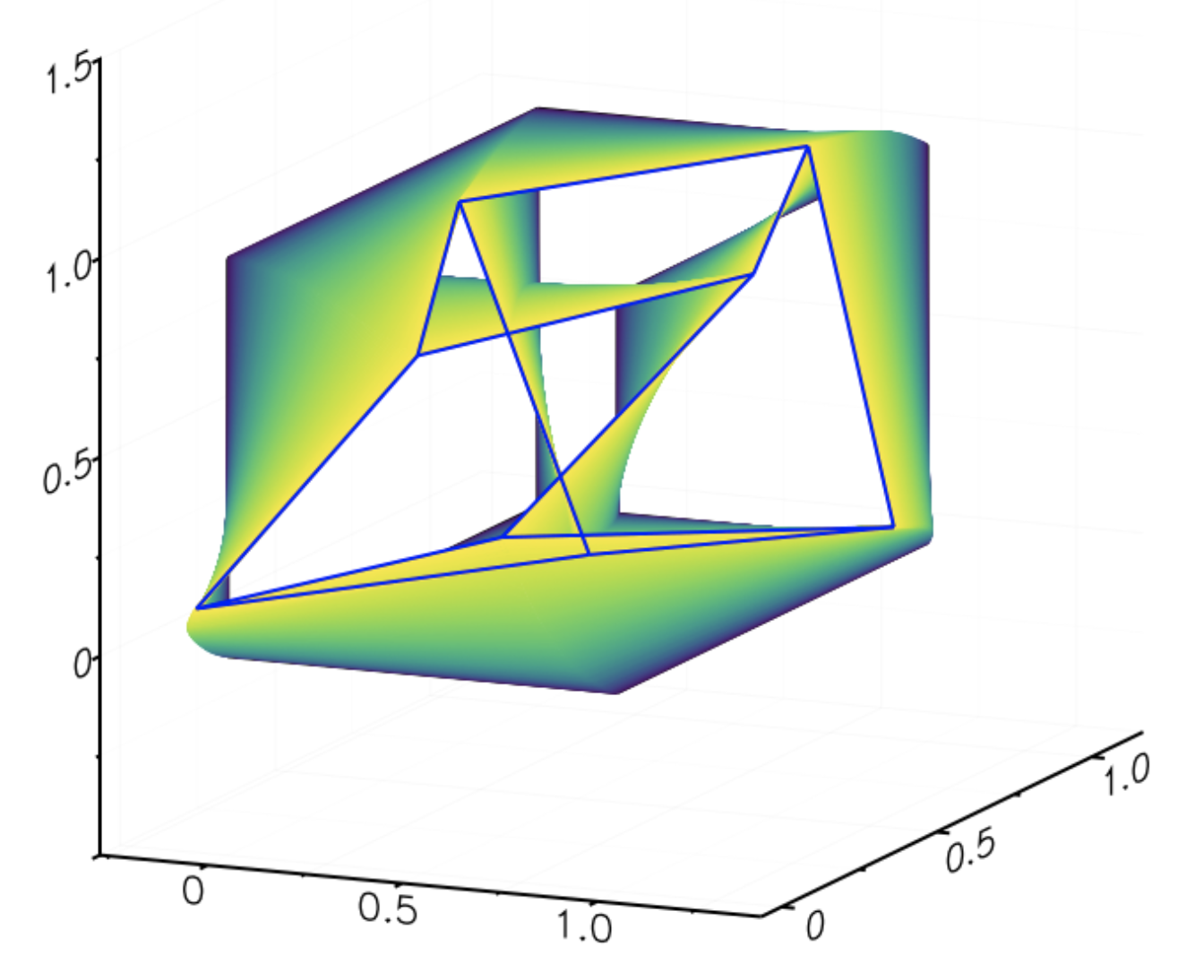}
  \caption{Cube flexing, from \cite{FrohmaderHeaton2020epsilonlocalrigidity}}\label{figure:flexible-cube}
\endminipage
\end{figure}

Recall in Section \ref{section:preliminary-setup} we evaluated the rank of the Jacobian $dg$ at randomly chosen values of $x_{ik}$, always finding it to be rank $12$ for the 3-prism. Since it has $18$ columns, we say it has \textit{corank} $6$. This means there is a $6$-dimensional nullspace. Since we know the $\binom{d+1}{2}$-dimensional group of rigid motions acts, we can immediately construct $6$ independent null vectors for any configuration $p$. The $3$ translations we leave as an exercise, but for the infinitesimal rotations, we apply a single skew-symmetric matrix to the coordinates of each node separately.
\begin{equation*}
    \left[ \begin{bmatrix} 0&-1&0\\1&0&0\\0&0&0 \end{bmatrix} \begin{bmatrix} x_{11}\\x_{12}\\x_{13} \end{bmatrix}, \begin{bmatrix} 0&-1&0\\1&0&0\\0&0&0 \end{bmatrix} \begin{bmatrix} x_{21}\\x_{22}\\x_{23} \end{bmatrix}, \, \, \dots \right]
\end{equation*}
There are three independent skew-symmetric matrices to use, and therefore three independent null vectors. If these are the only independent null vectors, we say that $p$ is \textit{infinitesimally rigid}. The random embeddings shown in Figure \ref{figure:random} are all infinitesimally rigid, while the embedding (\ref{equation:3prism-embedding-p0}) is not infinitesimally rigid, since there is a $7$th independent null vector, which (to three digits) is
\begin{multline}\label{equation:3-prism-flex-null-vector}
    \left(0.000,\,1.58,\,0.263,\,-1.37,\,-0.789,\,0.263,\,1.37,\,-0.789, \right. \\ \left. \,0.263,\,-0.789,\,1.37,\,-0.263,\,-0.789,\,-1.37,\,-0.263,\,1.58,\,0.000,\,-0.263\right) .
\end{multline}
Checking infinitesimal rigidity is a linearized version of the following. We say a tensegrity framework $p$ is \textit{rigid} if the only continuous deformations $p(t):[0,1] \to \mathbb{R}^{nd}$ preserving the member constraints (\ref{equation:member-constraints}) are rigid motions. Here $[0,1]$ is an interval of real numbers, representing time. Intuitively, $p(t)$ gives a movie of the framework deforming, or changing shape. If that deformation is equivalent to rotating or translating our reference frame, it is not really a change of shape, since no pairwise distance between nodes changed. We require $p(t)$ to preserve pairwise distances between nodes \textit{connected by an edge}, but if some non-edge distance can change, we say $p$ is \textit{flexible}. Thus a configuration is either rigid or flexible. Figures \ref{figure:flexible-honeycomb} and \ref{figure:flexible-cube} show flexible structures.

A configuration can be rigid without being infinitesimally rigid (linearization doesn't capture everthing). But, it turns out that if $p$ is infinitesimally rigid, then it is also rigid. To see why this is true, consider that $p$ is rigid if the \textit{local dimension} $\text{dim}_p V_{\mathbb{R}}(g)$ of the real algebraic set $V_{\mathbb{R}}(g)$ at the point $p$ is $\binom{d+1}{2}$. For precise definitions see \cite{WamplerHauensteinSommese2011mechanismMobility}. Since we always have rigid motions, we know that $\binom{d+1}{2} \leq \text{dim}_p V_{\mathbb{R}}(g)$. If we evaluate the rank (and hence corank) of $dg$ at the point $p$ and obtain $\text{corank } dg|_p = \binom{d+1}{2}$, then we know that the local complex dimension $\text{dim}_p V(g) \leq \binom{d+1}{2}$. But the local real dimension is always less than or equal to the local complex dimension, and since our polynomials have real coefficients they are equal at smooth points. So in fact we have shown that if $\text{corank }dg|_p = \binom{d+1}{2}$ then $p$ is a smooth point of $V(g)$ and also that $\text{dim}_p V_{\mathbb{R}}(g) = \binom{d+1}{2}$, so that $p$ is rigid. In conclusion, the \textit{singular locus} of $V(g)$ are the configurations which are \textit{interesting}, in the sense that infinitesimal rigidity tests may fail. We will return to this in Section \ref{section:adjacent-minors}.

The next degree after linear is quadratic. Every pure quadratic polynomial function $\mathbb{R}^N \to \mathbb{R}$ can be written as $x^T K x$ for some symmetric matrix $K$. If $x^T K x$ represents an energy function, then we (and nature) may want to find the minimum energy, if it exists. Here factorization is a powerful tool. Every matrix $K$ constructed as $K = B^T B$ for some rectangular matrix $B$ is positive semidefinite, meaning that $x^T K x \geq 0$ for all vectors $x \in \mathbb{R}^N$. To see why, write
\begin{equation*}
    x^T K x = x^T B^T B x = (Bx)^T (Bx) = \left\lVert Bx \right\rVert^2
\end{equation*}
which is nonnegative since the squared length of the vector $Bx$ will never be negative. Inserting a diagonal matrix whose entries are nonnegative preserves this property, so that if $K = B^T \text{diag}(c) B$ then $K$ is still positive semidefinite. Let $c = (c_{ij}) \in \mathbb{R}^m$ be a vector with $m$ components, one for every edge $ij \in E$. The stiffness matrix
\begin{equation}\label{equation:stiffness-matrix-Kc}
    K_c = (dg|_p)^T \text{diag}(c) (dg|_p)
\end{equation}
comes in exactly this form, so that with nonnegative $c_{ij}$ material constants, the function $x^T K_c x$ is positive semidefinite. In fact, the vectors giving zero $v^T K_c v = 0$ are the infinitesimal flexes and rigid motions $v$.

Definition 3.3.1 of \cite{ConnellyWhiteley1996secondOrderRigidityAndPrestressStabilityTensegrityFrameworks} gives a quadratic energy function $H_{w,c}$ that depends on a choice of vectors $w = (w_{ij}) \in \mathbb{R}^m$ and $c = (c_{ij}) \in \mathbb{R}^m$ as in
\begin{align}\label{equation:prestress-energy-function-Hwc}
    H_{w,c}(x) &= x^T\Omega_w x + x^T K_c x\\
     &= x^T \left( \Omega_w + K_c \right) x,\nonumber
\end{align}
where $K_c$ is from (\ref{equation:stiffness-matrix-Kc}) and $\Omega_w$ is the \textit{Kronecker product} of a weighted graph Laplacian with the identity matrix
\begin{equation}\label{equation:stress-matrix-Omegaw}
    \Omega_w = A^T \text{diag}(w) A \otimes I_d.
\end{equation}
Recall from Section \ref{section:understanding-rigidity-matrix} that the weighted graph Laplacian is constructed by using the incidence matrix $A$ of the graph, inserting a diagonal matrix of weights. For the 3-prism we have $\text{diag}(w)=$
\begin{equation*}\tiny
    \left[\begin{array}{rrrrrrrrrrrr}
w_{12} & 0 & 0 & 0 & 0 & 0 & 0 & 0 & 0 & 0 & 0 & 0 \\
0 & w_{13} & 0 & 0 & 0 & 0 & 0 & 0 & 0 & 0 & 0 & 0 \\
0 & 0 & w_{14} & 0 & 0 & 0 & 0 & 0 & 0 & 0 & 0 & 0 \\
0 & 0 & 0 & w_{15} & 0 & 0 & 0 & 0 & 0 & 0 & 0 & 0 \\
0 & 0 & 0 & 0 & w_{23} & 0 & 0 & 0 & 0 & 0 & 0 & 0 \\
0 & 0 & 0 & 0 & 0 & w_{25} & 0 & 0 & 0 & 0 & 0 & 0 \\
0 & 0 & 0 & 0 & 0 & 0 & w_{26} & 0 & 0 & 0 & 0 & 0 \\
0 & 0 & 0 & 0 & 0 & 0 & 0 & w_{34} & 0 & 0 & 0 & 0 \\
0 & 0 & 0 & 0 & 0 & 0 & 0 & 0 & w_{36} & 0 & 0 & 0 \\
0 & 0 & 0 & 0 & 0 & 0 & 0 & 0 & 0 & w_{45} & 0 & 0 \\
0 & 0 & 0 & 0 & 0 & 0 & 0 & 0 & 0 & 0 & w_{46} & 0 \\
0 & 0 & 0 & 0 & 0 & 0 & 0 & 0 & 0 & 0 & 0 & w_{56}
\end{array}\right], \, \, A = \left[\begin{array}{rrrrrr}
1 & -1 & 0 & 0 & 0 & 0 \\
1 & 0 & -1 & 0 & 0 & 0 \\
1 & 0 & 0 & -1 & 0 & 0 \\
1 & 0 & 0 & 0 & -1 & 0 \\
0 & 1 & -1 & 0 & 0 & 0 \\
0 & 1 & 0 & 0 & -1 & 0 \\
0 & 1 & 0 & 0 & 0 & -1 \\
0 & 0 & 1 & -1 & 0 & 0 \\
0 & 0 & 1 & 0 & 0 & -1 \\
0 & 0 & 0 & 1 & -1 & 0 \\
0 & 0 & 0 & 1 & 0 & -1 \\
0 & 0 & 0 & 0 & 1 & -1
\end{array}\right].
\end{equation*}
Then $A^T \text{diag}(w) A$ becomes
\begin{equation*} \footnotesize
    \left[\begin{array}{rrrrrr}
\begin{array}{c}
    w_{12} + w_{13} \\+ w_{14} + w_{15}
\end{array} & -w_{12} & -w_{13} & -w_{14} & -w_{15} & 0 \\
-w_{12} & \begin{array}{c}
    w_{12} + w_{23} \\+ w_{25} + w_{26}
\end{array} & -w_{23} & 0 & -w_{25} & -w_{26} \\
-w_{13} & -w_{23} & \begin{array}{c}
    w_{13} + w_{23} \\+ w_{34} + w_{36}
\end{array} & -w_{34} & 0 & -w_{36} \\
-w_{14} & 0 & -w_{34} & \begin{array}{c}
    w_{14} + w_{34} \\+ w_{45} + w_{46}
\end{array} & -w_{45} & -w_{46} \\
-w_{15} & -w_{25} & 0 & -w_{45} & \begin{array}{c}
    w_{15} + w_{25} \\+ w_{45} + w_{56}
\end{array} & -w_{56} \\
0 & -w_{26} & -w_{36} & -w_{46} & -w_{56} & \begin{array}{c}
    w_{26} + w_{36} \\+ w_{46} + w_{56}
\end{array}
\end{array}\right].
\end{equation*}
The Kronecker product with an identity matrix simply \textit{spreads out} the other matrix across all spatial dimensions. For a small example consider
\begin{equation*}
    \left[ \begin{array}{cc}
        1 & 2 \\
        3 & 4
    \end{array} \right] \otimes I_3 = \left[ \begin{array}{cccccc}
        1 & & &2 & &  \\
         &1& & &2 & \\
         & &1& & &2 \\
        3 & & &4 & &  \\
         &3& & &4 & \\
         & &3& & &4 \\
    \end{array} \right].
\end{equation*}
Definition 3.3.1 of \cite{ConnellyWhiteley1996secondOrderRigidityAndPrestressStabilityTensegrityFrameworks} says that a tensegrity framework is \textit{prestress rigid} if there is a self stress $w^T dg|_p = 0$ with all $w_{ij} \neq 0$ and a positive vector $c > 0$ such that the energy function $H_{w,c}$ of (\ref{equation:prestress-energy-function-Hwc}) is positive semidefinite and gives zero only on infinitesimal rigid motions. Intuitively, the energy function $H_{w,c}$ builds in stabilizing equilibrium tensions $w$. We interpret the weights $w_{ij}$ appearing in $\Omega_w$ as tensions on each edge. If $w^T dg|_p = 0$ then the internal forces (as in Section \ref{section:understanding-rigidity-matrix}) exactly balance at every node. Therefore, the configuration $p$ has no reason to move, at least not due to the tensions $w_{ij}$. Without the tensions $w$ the energy is $x^T K_c x$ only, and infinitesimal flexes allow movement with no energy increase. But with tensions $w$, since $H_{w,c}$ is positive definite on the flexes, now those flexes cost energy! Thus the tensions \textit{neutralize} the flexes.

In \cite{ConnellyWhiteley1996secondOrderRigidityAndPrestressStabilityTensegrityFrameworks} they prove that prestress rigidity implies rigidity. They also show we can ignore $K_c$ and instead check a simpler condition involving $\Omega_w$ alone. Proposition 3.4.2 of \cite{ConnellyWhiteley1996secondOrderRigidityAndPrestressStabilityTensegrityFrameworks} shows that if there exists some $w^T dg|_p = 0$ such that $\Omega_w$ is positive definite on any flex subspace $F$ with $\text{Null }dg|_p = R \oplus F$ from Equation (\ref{equation:nullspace-decomposition-R-plus-F}), then the tensegrity framework $p$ is prestress rigid. This turns the problem of checking prestress rigidity of a framework into a semidefinite program (SDP). For an introduction to SDPs, see Chapter 12 of \cite{MS2019}, where they explain that an SDP is the non-abelian version of a linear program (LP). Specifically, LPs are a subset of SDPs because diagonal matrices with positive entries are a subset of square matrices with positive eigenvalues.

Let $F$ denote a rectangular matrix whose columns span the infinitesimal flexes of $\text{Null }dg|_p = R \oplus F$, abusing notation. Assume we have a framework with a $2$-dimensional space of $w^T dg|_p = 0$ (analogous cases for higher dimensions will be clear). Assume also that vectors $w_1, w_2 \in \mathbb{R}^m$ form a basis. Then we solve an SDP to find scalars $a_1, a_2 \in \mathbb{R}$ such that
\begin{equation}\label{equation:prestress-SDP}
    F^T \left( a_1 \cdot \Omega_{w_1} + a_2 \cdot \Omega_{w_2} \right) F \succ 0
\end{equation}
is a positive definite matrix. This means our configuration $p$ is a local minimum (modulo rigid motions) for some energy function $H_{w,c}$. If the search succeeds then $a_1 w_1 + a_2 w_2$ gives a self stress witnessing the prestress rigidity of $p$. Such a semidefinite program is solvable in polynomial time by efficient interior point algorithms, making prestress rigidity an attractive criterion for deciding rigidity of $p$. This SDP is even more attractive when there is only one basis vector $w_1$, and only one infinitesimal flex $v_1$. Then (\ref{equation:prestress-SDP}) becomes
\begin{equation*}
    v_1^T \left( a_1 \cdot \Omega_{w_1} \right) v_1 > 0.
\end{equation*}
We can solve the SDP almost by hand, checking $a_1 \in \{1,-1\}$ and in either case computing the single scalar output and observing whether it is positive or not.

For the 3-prism this is the case; the left nullspace $w^T dg|_p = 0$ is only 1-dimensional. To three digits, we show the entries of a self stress vector $w$ in a table so that each entry comes with its associated edge label.
\begin{equation}\label{equation:table-of-prestress-entries-and-edges}
    \begin{tabular}{|c|c|c|c|c|c|c|c|c|c|c|c|} \hline
$\left(1, 2\right)$ & $\left(1, 3\right)$ & $\left(1, 4\right)$ & $\left(1, 5\right)$ & $\left(2, 3\right)$ & $\left(2, 5\right)$ & $\left(2, 6\right)$ & $\left(3, 4\right)$ & $\left(3, 6\right)$ & $\left(4, 5\right)$ & $\left(4, 6\right)$ & $\left(5, 6\right)$ \\ \hline
$1.00$ & $1.00$ & $-1.73$ & $1.73$ & $1.00$ & $-1.73$ & $1.73$ & $1.73$ & $-1.73$ & $1.00$ & $1.00$ & $1.00$ \\ \hline
\end{tabular}
\end{equation}
We also have a one-dimensional space of infinitesimal flexes, given by the span of the vector $v$ of (\ref{equation:3-prism-flex-null-vector}). The code provided at \cite{heaton-tensegrity-trusses-github} evaluates
\begin{equation*}
    v^T \Omega_{w} v = 89.56922 > 0
\end{equation*}
showing that $\Omega_w$ for the self stress (\ref{equation:table-of-prestress-entries-and-edges}) is positive definite on the space of infinitesimal flexes $F = \text{span}\{v\}$. We have solved a very small semidefinite program. Since we have a tensegrity framework, we must also check that cables have positive tensions, and struts have negative tensions, while bars can have either (see the edge labels in Figure \ref{figure:3prism-flexes}). We conclude the 3-prism is prestress rigid.

\section{Numerical algebraic geometry}\label{section:numerical-algebraic-geometry}

An interesting task is to follow the motion of planets through space by solving a system of ordinary differential equations (ODEs) in a variable $t$. We must keep track of the $x,y,z$ coordinates of the planet starting from its \textit{initial conditions}, the $x,y,z$ coordinates at $t=0$. In the same way, consider the polynomial $g = x^3 - 1$, and picture its three solutions as the initial conditions of three planets in the plane (our own version of the novel by Liu Cixin). Here, we are identifying $\mathbb{C}$ with $\mathbb{R}^2$. We could imagine following the motion of these three solutions just as the Trisolarans followed their three suns. But what makes solutions move, if not gravity?

\begin{figure}[!htb]
    \centering
    \includegraphics[width=0.4\textwidth]{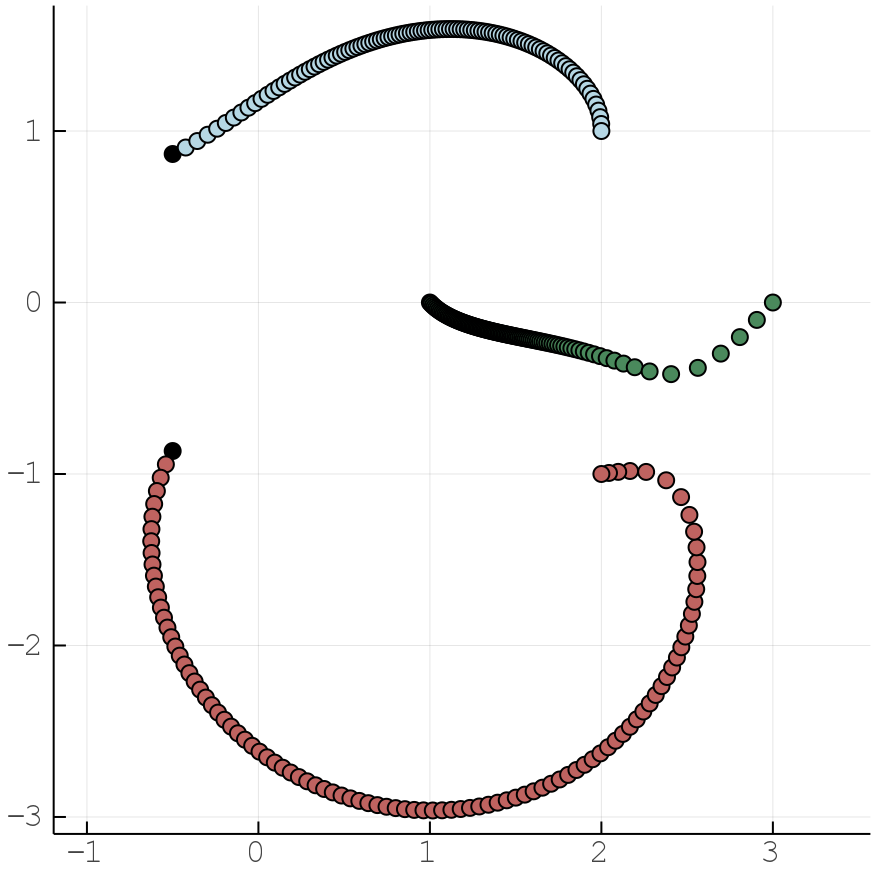}
    \includegraphics[width=0.4\textwidth]{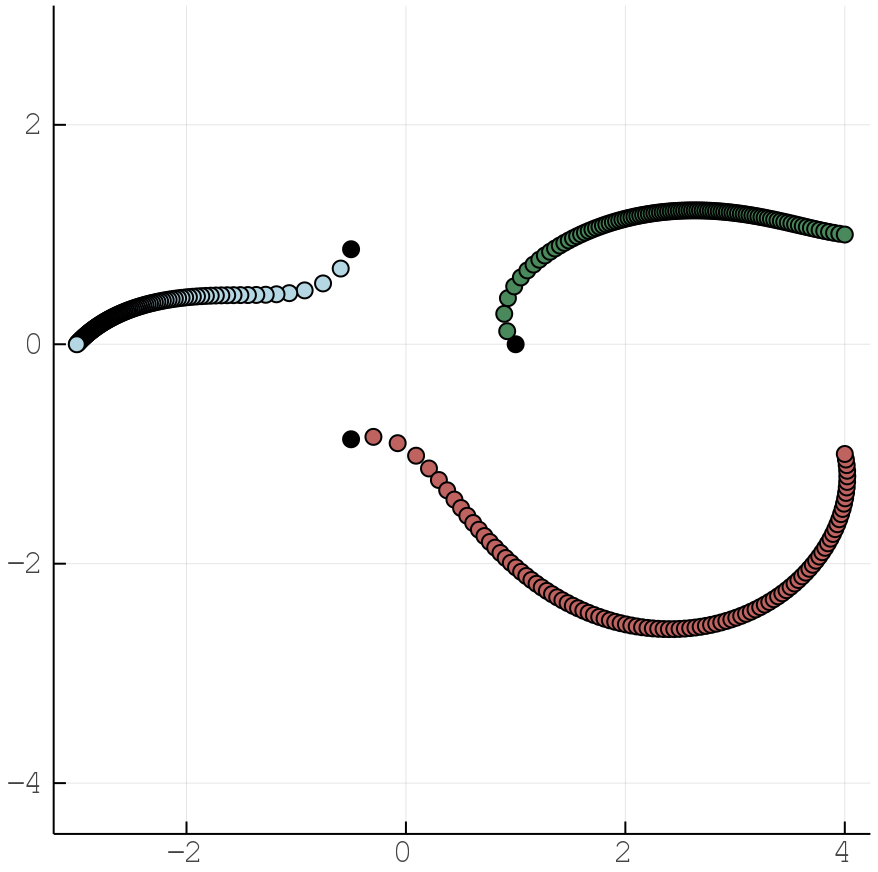}
    \caption{Following solution paths like following planets in space. Computed in \texttt{Julia} \cite{Julia2017}}
    \label{figure:basic-homotopy-planets}
\end{figure}

Consider another polynomial $f = x^3 - 7x^2 + 17x - 15$, which, by the Fundamental Theorem of Algebra, has three roots (possibly with multiplicity). Imagine I didn't tell you that $f = (x-3)(x-(2+i))(x-(2-i))$, and we would like to discover these three roots. The subject of \textit{numerical algebraic geometry} \cite{BatesSommeseHauensteinWampler2013BertiniTEXT, HauensteinSommese2017whatisNumericalAlgebraicGeometry, SommeseWampler2005numericalSolutionofSystemsofPolynomialsTEXT} is concerned with solving this problem by perturbing known solutions (planets with initial conditions) towards the unknown. More generally, it develops and applies numerical algorithms to understand algebraic sets, drawing on the theory of algebraic geometry. Figure \ref{figure:basic-homotopy-planets} contains the end results of code available at \cite{heaton-tensegrity-trusses-github}. The example on the left shows the motion of the three unit circle solutions of $x^3 -1$ moving towards the zeros $3, 2+i, 2-i$. The example on the right shows the motion for another polynomial $f = x^3 -5x^2 -7x +51$, whose roots we discover as $-3, 4+i, 4-i$.

Given a start system $g(x) = x^3 - 1$ whose solutions we know, we can perturb it towards a target system, like $f(x) = x^3 - 7x^2 + 17x - 15$, whose solutions we want to discover. This is accomplished by using \textit{homotopy continuation}. We will not explain this fully, but give the basic idea. Consider
\begin{equation*}
    h(x,t) = (1-t)f + t g.
\end{equation*}
Numerical algebraic geometers believe time starts at $t=1$ and runs backwards toward $t=0$. But regardless of time's flow, at $t=1$ we have $h(x,1) = g(x)$, whose solutions we know. If we follow the motion of these ``planets'' towards $t=0$ we will discover solutions to $h(x,0) = f(x)$, which were previously unknown. It turns out that multiplying by a random point on the unit circle $\gamma \in S^1 \subset \mathbb{C}$ allows us to avoid \textit{singularities} in parameter space, so actually we use the homotopy which has been twisted by the so-called \textit{gamma trick}, as in
\begin{equation*}
    h(x,t) = (1-t)f + \gamma t g.
\end{equation*}
Taking $x(t)$ to be the unknown solution path (the trajectory of a planet) we have that $h(x(t),t) = 0$ for all $t$, then differentiating with respect to $t$ and applying the chain rule gives the \textit{Davidenko differential equation}
\begin{equation*}
    \frac{\partial h}{\partial x} \frac{dx}{dt} + \frac{\partial h}{\partial t} = 0.
\end{equation*}
For systems of equations in $N$ complex variables $\frac{\partial h}{\partial x}$ denotes the Jacobian of $h$ with respect to the variables $x_1, \dots, x_N$, while the vectors $\frac{dx}{dt}$ and $\frac{\partial h}{\partial t}$ have $N$ components each. For full details and wonderfully explained examples on numerical algebraic geometry we refer the reader to either of the textbooks \cite{BatesSommeseHauensteinWampler2013BertiniTEXT, SommeseWampler2005numericalSolutionofSystemsofPolynomialsTEXT} or also the \texttt{HomotopyContinuation.jl} website \cite{BreidingTimme2018homotopyContinuationJL}. This is a \texttt{Julia} implementation of homotopy continuation which we use for all our examples.

Using theoretical results from algebraic geometry as its foundation, numerical algebraic geometry aims to provide algorithms that produce \textit{all solutions} reliably and efficiently. Since we have systems of polynomial equations, Bezout's theorem gives an upper bound on the total possible number of solutions. Using appropriate starting systems will result in finding all solutions, in contrast to methods like Newton-Raphson. In addition, through the \textit{numerical irreducible decomposition}, we can understand the positive-dimensional solution components geometrically, something that other numerical solving techniques cannot achieve. In particular, the concept of \textit{witness sets} provides a finite data structure representing each irreducible component $X$. We keep track of the polynomials $f$, a linear space of complementary dimension $L$, and the \textit{degree(X)}-many witness points, intersections of $L$ with $V(f)$.

In the beginning, research on homotopy continuation focused on computing isolated solutions by \textit{total degree homotopy}. Theorems from algebraic geometry showed we can compute all isolated solutions, and also how to avoid singularities in parameter space. Bezout's theorem gives an upper bound on the total number of isolated solutions: the product of the degrees of all the polynomials involved. For a system we will solve in a moment, this number will be $67,108,864$. This is the number of planets whose paths we must track. Later, using more theoretical results from algebraic geometry, the algorithms improved even further. For example, a theorem of Shafarevich on multi-homogeneous systems can bring that number down to $8,503,056$, far fewer planets \cite{Shafarevich2013BasicAlgebraicGeometryTEXT}. Bernstein's theorem connects the number of solutions to the mixed volume of \textit{Newton polytopes} \cite{Bernstein1975numberRootsOfASystemOfEquations}. This was used by Huber and Sturmfels to construct \textit{polyhedral homotopies} \cite{HuberSturmfels1995PolyhedralMethodForSolvingSparsePolynomialSystems}. In our example this is the winner, yielding only $1,062,880$ paths to track. The recent paper \cite{BreidingSturmfelsTimme20203264conicsPerSeconds} discusses how numerical algebraic geometry can be applied to \textit{Steiner's problem}, and in particular finds all $3264$ conics tangent to a given five conics, in just one second. We recommend this paper for a demonstration of how numerical methods and enumerative geometry complement each other.

To use numerical algebraic geometry to compute configurations of a bar framework with $V(g) \subset \mathbb{C}^{nd}$ we first realize that we would like to compute modulo the group action. We draw ispiration from the method of \textit{moving frames} \cite{OlverMovingFrames}, which often uses \textit{coordinate cross-sections}. For us, we realize that given any configuration $(x_{ik}) \in \mathbb{R}^{nd}$ we can \textit{change coordinates} on $\mathbb{R}^d$ so that in the new coordinates, we have more zeros. First, we can always translate the origin to node 1. Then, we can rotate coordinates until node 2 is on the $x$-axis. Finally, we can further rotate about the $x$-axis until node 3 lies in the $(x,y)$-plane. For the configuration (\ref{equation:3prism-embedding-p0}) of the 3-prism, we obtain
\begin{equation*}
\left[ \begin{array}{ccc}
        1 & 0 & 0\\ -\frac{1}{2} & \frac{\sqrt{3}}{2} & 0\\  -\frac{1}{2} & -\frac{\sqrt{3}}{2} & 0\\ -\frac{\sqrt{3}}{2} & -\frac{1}{2} & 3\\ \frac{\sqrt{3}}{2} & -\frac{1}{2} & 3\\ 0 & 1 & 3
    \end{array} \right] \mapsto 
\left[
\begin{array}{ccc}
0.0 & 0.0 & 0.0 \\
1.7320508075688772 & 0.0 & 0.0 \\
0.8660254037844388 & -1.5 & 0.0 \\
1.3660254037844386 & -1.3660254037844386 & 3.0 \\
-0.1339745962155613 & -0.5 & 3.0 \\
1.3660254037844388 & 0.3660254037844386 & 3.0 \\
\end{array}
\right].
\end{equation*}
This change of coordinates can be computed easily using the $QR$ matrix factorization, as we do in \cite{heaton-tensegrity-trusses-github}. In this way, we can reduce from $nd$ variables to only $N = nd - \binom{d+1}{2}$ variables. The only difference from Figure \ref{figure:basic-homotopy-planets} is that now we are following planets in real $2\cdot N$ dimensional space.

\begin{figure}[!htb]
    \centering
    \includegraphics[width=0.47\textwidth]{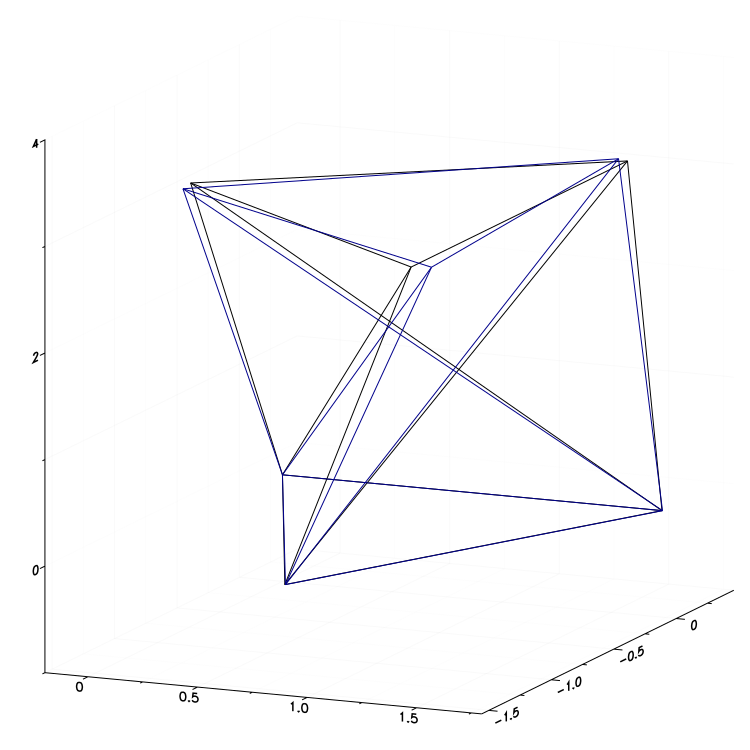}
    \includegraphics[width=0.47\textwidth]{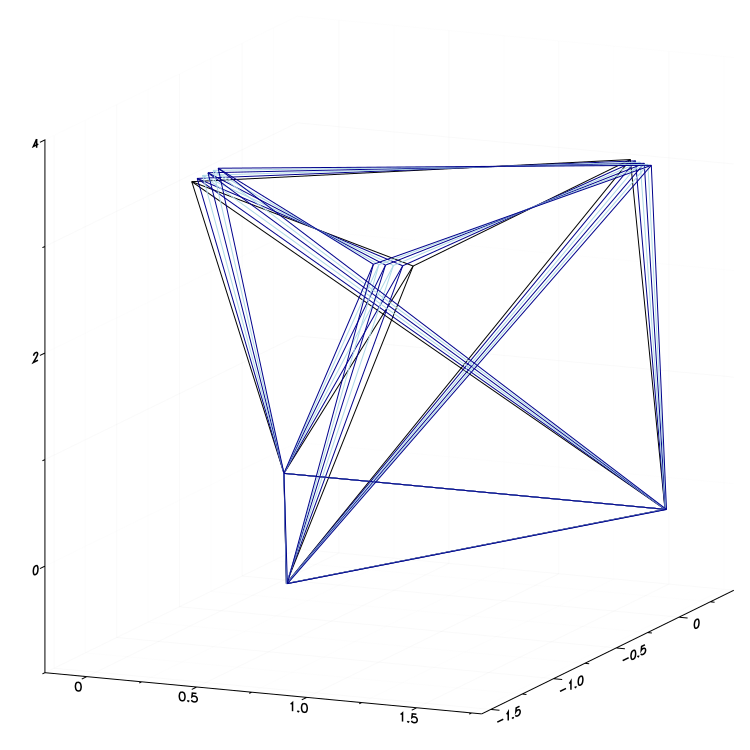}
    \caption{The 3-prism deforms}
    \label{figure:3prism-deforms}
\end{figure}

Glance at Figure \ref{figure:3prism-deforms}.
Using \texttt{HomotopyContinuation.jl} \cite{BreidingTimme2018homotopyContinuationJL} we can compute solutions $p(0+{\Delta t})$ which are very close to our initial configuration $p(0) = p$ in (\ref{equation:3prism-embedding-p0}). Here is one of the first deformations where the 3-prism begins to twist downwards, as you can see in the $z$-coordinates of nodes $4,5,6$ in (\ref{equation:downwards}).
\begin{equation}\label{equation:downwards}
\left[
\begin{array}{ccc}
0.0 & 0.0 & 0.0 \\
1.7345015098619578 & 0.0 & 0.0 \\
0.868440136662386 & -1.4992275136004456 & 0.0 \\
1.434394418877309 & -1.322820159782012 & 2.9867578031247515 \\
-0.12780675639331027 & -0.5722571170386941 & 2.9884237516052568 \\
1.3032879403929745 & 0.40433597324593706 & 2.9865056925919635 \\
\end{array}
\right]
\end{equation}
The 3-prism can also untwist upwards, as seen in (\ref{equation:upwards}).
\begin{equation}\label{equation:upwards}
\left[
\begin{array}{ccc}
0.0 & 0.0 & 0.0 \\
1.7366579715554198 & 0.0 & 0.0 \\
0.8689191994267838 & -1.5000751141689224 & 0.0 \\
1.2441020030189796 & -1.4251233353656827 & 3.022965523609326 \\
-0.12034129224092607 & -0.35340935623990566 & 3.024004212778408 \\
1.4887072950829638 & 0.2911265336721315 & 3.021803451042337 \\
\end{array}
\right]
\end{equation}

How did we compute (\ref{equation:downwards}) and (\ref{equation:upwards})? We start with the system of equations $g(x) = 0$ in (\ref{equation:member-constraints}) but then adjoin one additional equation:
\begin{equation*}
    \ell_v = v^T x - v^T p = 0.
\end{equation*}
Here, we can choose $v \in \mathbb{R}^N$ randomly, or we could choose $v$ from some infinitesimal flex. Either way, the equation $\ell_v(x) = 0$ is solved by the initial configuration $p \in \mathbb{R}^N$. We then \textit{perturb} this equation to
\begin{equation*}
    \ell_{v,\epsilon} = v^T x - v^T p - \epsilon = 0.
\end{equation*}
for some small real $\epsilon \in \mathbb{R}$. Geometrically, we start with a hyperplane passing through the point $p$ and then slightly move the hyperplane in the direction of its normal vector $v$. If $v$ is an infinitesimal flex, we are moving the hyperplane in the direction of that flex in $\mathbb{R}^N$. In the computer, we use a \textit{real parameter homotopy} (without $\gamma$) as in
\begin{equation*}
    h_t(x) = (1-t)\left[ \begin{array}{c}
        g \\
        \ell_{v,\epsilon}  
    \end{array} \right] + t \left[ \begin{array}{c}
        g \\
        \ell_v  
    \end{array} \right]
\end{equation*}
to follow our solution $p$ (only one planet to follow) from $t=1$ to $t=0$. Here, $g$ represents all our member constraints $g_{ij}$ for $ij \in B$. Repeating this many times, we can compute points like those in (\ref{equation:downwards}) and (\ref{equation:upwards}).

We had concluded the 3-prism was prestress rigid in Section \ref{section:rigidity-prestress-stability}. Why is it deforming in Figure \ref{figure:3prism-deforms}? Are these points (\ref{equation:downwards}) and (\ref{equation:upwards}) really solutions to $g(x) = 0$? In the reality of numerical computation, we give up on \textit{exact solutions}. Instead, we make reasonable decisions based on the results of our calculations. In particular, using the code at \cite{heaton-tensegrity-trusses-github} you can attempt to compute these deformations of the 3-prism yourself using \texttt{HomotopyContinuation.jl} \cite{BreidingTimme2018homotopyContinuationJL}. You will find that most often the computation fails. It tells you, ``Sorry, there are no real-valued solutions to your parameter homotopy.'' But if you are \textit{extremely persistent}, repeatedly choosing new hyperplanes and trying again, you will eventually find points like (\ref{equation:downwards}) and (\ref{equation:upwards}). We argue that although these points are not real-valued solutions to $g(x)=0$, they do have meaning! First of all, they correspond to complex-valued solutions with very small imaginary parts. It is no accident that the nodes of these newly computed configurations have been perturbed in the same direction as the infinitesimal flexes of Figure \ref{figure:3prism-flexes}. In fact, if you build the 3-prism yourself, you will find that it deforms a small amount in exactly that way. No bar is perfectly rigid. Therefore, even though (\ref{equation:downwards}) and (\ref{equation:upwards}) are not exactly points on $V_{\mathbb{R}}(g)$, they are still very useful to compute.

\begin{figure}[!htb]
    \centering
    \includegraphics[width=0.9\textwidth]{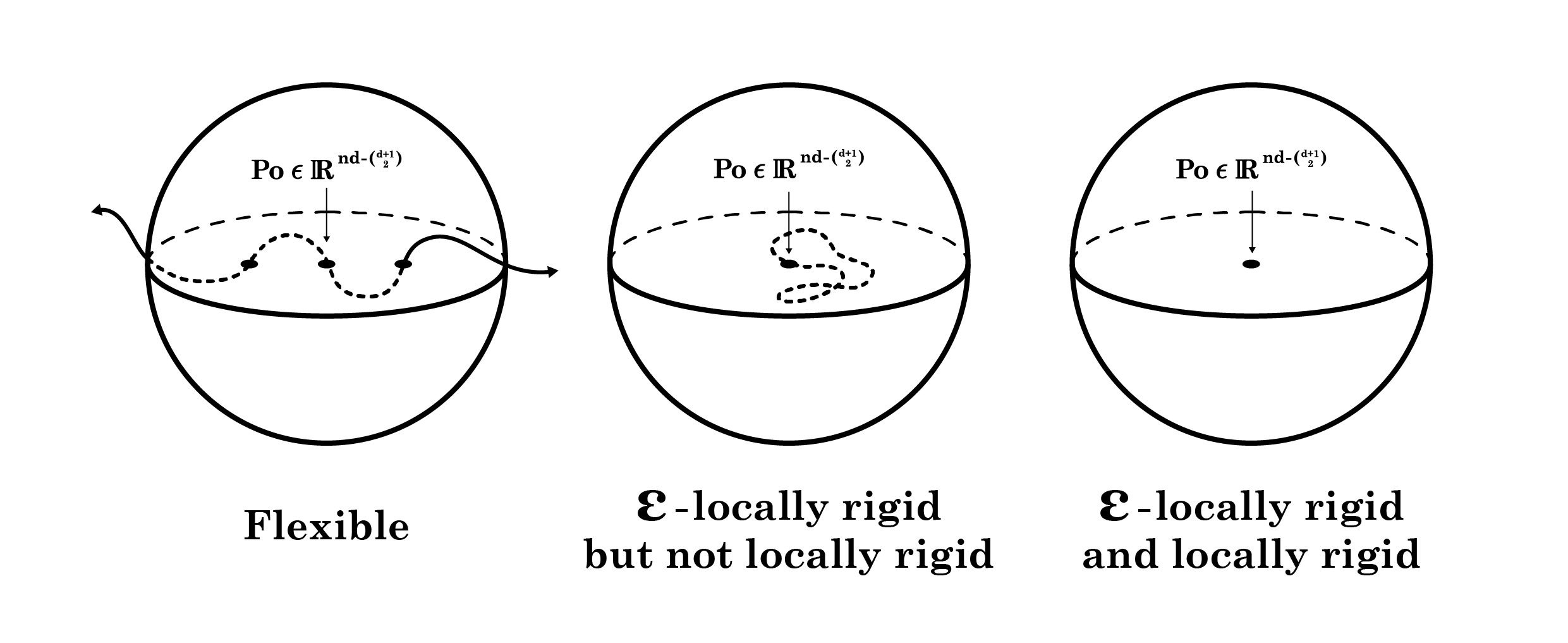}
    \caption{Epsilon local rigidity}
    \label{figure:epsilon-local-rigidity}
\end{figure}

Say you observe small deformations in some other example, but you still suspect your configuration $p$ is locally rigid. In that case, we can use \textit{epsilon local rigidity} of \cite{FrohmaderHeaton2020epsilonlocalrigidity}. In this calculation, we add one more equation to the system $g$, namely the equation $s_\varepsilon(x)=0$ of an $\varepsilon$-sphere centered at $p$ for some $0 < \varepsilon \in \mathbb{R}$ of your choosing. Consider Figure \ref{figure:epsilon-local-rigidity}.
If $p$ is locally rigid then small $\varepsilon$-spheres should not intersect $V_{\mathbb{R}}(g)$ at all. On the other hand, if there is a deformation $p(t)$, then there should be real-valued intersection points between the sphere and $V_{\mathbb{R}}(g)$.

Theorem 7 and Algorithm 1 of \cite{FrohmaderHeaton2020epsilonlocalrigidity} provide theoretical guarantees for computing these real-valued intersection points. If we find no real-valued solutions to $V_\mathbb{R}(g) \cap V_\mathbb{R}(s_\varepsilon)$, we have shown that $p$ is $\varepsilon$-locally rigid. The computation sums the squares of all the member constraints, obtaining a pure $(N-1)$-dimensional hypersurface. It then uses \textit{Lagrange multipliers} to minimize Euclidean distance from some randomly chosen point. Results from \cite{AubryRouillierSafeyElDin2002, Hauenstein2012, RouillierRoySafeyElDin2000, Seidenberg1954} show that minimizing Euclidean distance guarantees at least one point on each connected component of the real algebraic set. Therefore, if we find none, we can conclude that the $\varepsilon$-sphere misses $V_{\mathbb{R}}(g)$ entirely. In that case all deformations $p(t)$, should they exist, are contained in the $\varepsilon$-ball. In applications, deformations that stay so close to the initial configuration can probably be ignored. The study of \textit{real algebraic geometry} is extremely interesting and also difficult, as many of the nice results from algebraic geometry over $\mathbb{C}$ can drastically fail for $\mathbb{R}$.

We made this calculation for the 3-prism, obtaining $\varepsilon$-local rigidity for $\varepsilon=0.1$, even though the deformations displayed in Figure \ref{figure:3prism-deforms} were outside such a small $\varepsilon$-ball. This confirms what we already knew from prestress rigidity, that the deformations of Figure \ref{figure:3prism-deforms} do not exactly satisfy the member constraints. As we mentioned earlier, this example demonstrates the usefulness of algebraic geometry for numerical calculations, just like the theory of linear algebra informs numerical linear algebra. Using a naive total degree homotopy would require tracking $67,108,864$ paths. Using a theorem of Shafarevich brings that number down to $8,503,056$. Using mixed volumes of Newton polytopes lowers it further: $1,062,880$. Tracking this last number of paths took just a few hours on my personal computer. You can try as well \cite{heaton-tensegrity-trusses-github}.

\section{Gr\"obner bases and primary decomposition}\label{section:adjacent-minors}

After discussing numerical, floating point calculations in Section \ref{section:numerical-algebraic-geometry}, we now consider a technique available with \textit{exact symbolic computation}, and in particular, Gr\"obner bases. As is the pattern, we will not provide all details. Briefly, Gr\"obner bases generalize three important algorithms: (1) \textit{Gaussian elimination} in linear algebra, (2) the \textit{Euclidean algorithm} for finding the greatest common divisor, and (3) the \textit{simplex method} of linear programming. The power of this approach is best appreciated in an example, so we give two. First, the particularly beautiful example of \textit{adjacent minors}: Recall that the configurations $(x_{ik})$ causing the Jacobian matrix to drop rank are especially interesting. The configuration (\ref{equation:3prism-embedding-p0}) of the 3-prism is such an example. These rank-dropping conditions are expressed by the vanishing of minors of matrices. Consider an $m \times n$ matrix of variables $x_{ij}$ as in
\begin{equation}\label{equation:2by5matrix}
    \left[ \begin{array}{ccccc}
        x_{11} & x_{12} & x_{13} & x_{14} & x_{15} \\
        x_{21} & x_{22} & x_{23} & x_{24} & x_{25} 
    \end{array} \right]
\end{equation}
for the $2 \times 5$ case. The variety of all matrices of rank $< r$ is an irreducible variety whose prime ideal is generated by all the $r \times r$ minors of the $m \times n$ matrix $x_{ij}$. One way to get reducible varieties determined by minors of matrices is to restrict to \textit{structured matrices} with a sparsity pattern determined by a graph $([n],E)$, as is the topic of this paper. These varieties will be more complicated, having several irreducible components of different dimensions and degrees. Another way to get reducible varieties from minors of matrices is to restrict our attention to only the \textit{adjacent minors}. These examples find application in statistics, and were considered in \cite{DiaconisEisenbudSturmfels1998latticeWalksPrimaryDecomposition, HostenShapiro2000primaryDecompositionLatticeBasisIdeals}, for example. Here we use them to introduce an algebraic technique called \textit{primary decomposition}, which will help us make sense of polynomial equations.

\begin{figure}[h]
    \centering
    \includegraphics[width=0.7\textwidth]{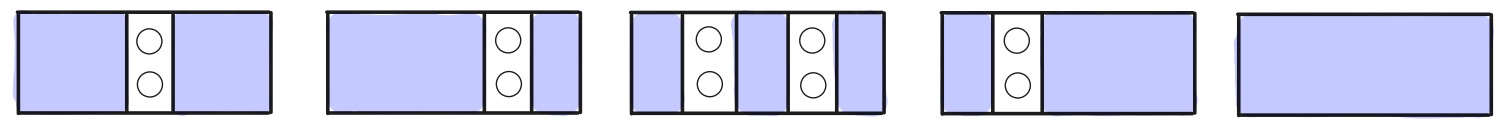}
    \caption{Irreducible components}
    \label{figure:adjacent-ranks}
\end{figure}

An ideal $I \subset \mathbb{Q}[x_1,\dots,x_n]$ is first of all a subset. Second, it is a vector subspace. Third, it is a module closed under linear combinations with arbitrary polynomial coefficients. If polynomials $g_1$ and $g_2$ are in $I$, then so is $f_1g_1 + f_2 g_2$, for any arbitrary polynomials $f_1,f_2 \in \mathbb{Q}[x_1,\dots,x_n]$. Every ideal $I \subset \mathbb{Q}[x_1,\dots,x_n]$ admits a primary decomposition
\begin{equation*}
    I = Q_1 \cap Q_2 \cap \cdots \cap Q_r
\end{equation*}
for primary ideals $Q_i$ with each $Rad Q_i$ a distinct \textit{associated prime}. For precise definitions, we refer the reader to \cite{CoxLittleOshea2015IdealsVarietiesAlgorithmsTEXT} or Chapter 3 of \cite{MS2019}. For our purposes, it is enough to say that breaking up our algebraic variety into irreducible pieces corresponds to ``breaking up'' our ideal as an intersection of primary ideals. An example will illustrate best. In the $2 \times 5$ matrix (\ref{equation:2by5matrix}) above, the ideal generated by all $2 \times 2$ minors is prime and its irreducible variety consists of all $2 \times 5$ matrices of rank $< 2$. If instead we only require that \textit{adjacent minors} vanish, we get a bigger variety with several irreducible components. The \texttt{SAGE} code is so short we can show it here:
\begin{verbatim}
xvarz = [var('x%s%s'%(i,j)) for i in range(1,2+1) for j in range(1,5+1)]
adjacentminors = [x11*x22 - x12*x21, x12*x23 - x13*x22,
                x13*x24 - x14*x23, x14*x25 - x15*x24]
R = PolynomialRing(QQ,xvarz)
I = R.ideal(adjacentminorz)
for P in I.associated_primes():
    print P; print;
\end{verbatim}
Introducing a shorthand for minors as in $13 := x_{11}x_{23} - x_{13}x_{21}$, we can describe the output of the code above by
\begin{equation*}
    \begin{array}{c}
        P_1 = \langle 12, \, x_{13}, \, x_{23}, \, 45 \rangle \\
        P_2 = \langle 12, \, 13,\, 23,\, x_{14}, \, x_{24} \rangle \\
        P_3 = \langle x_{12}, \, x_{22},\, x_{14},\, x_{24} \rangle \\
        P_4 = \langle x_{12}, \, x_{22},\, 34,\, 35,\, 45 \rangle \\
        P_5 = \langle 12, \, 13, \, 14, \, 15, \, 23, \, 24, \, 25, \, 34, \, 35,\, 45 \rangle
    \end{array}
\end{equation*}
Since we took $2 \times 2$ minors, the ideals are binomial ideals, and so primary decomposition is well-behaved (see \cite{EisenbudSturmfels1996binomialIdeals} for more information). In particular all the associated primes $P_1,\dots, P_5$ are also binomial ideals. In fact, we can understand them now. $P_1$ tells us the third column is zero, while the first and last $2 \times 2$ blocks are rank 1 matrices. We leave it as an exercise to describe all matrices satisfying each of the other ideals $P_2,P_3,P_4,P_5$, but leave a hint in Figure \ref{figure:adjacent-ranks}.
We also quote from \cite{Sturmfels2002solvingsystemsofpolynomialequationsTEXT}: ``If all adjacent $2 \times 2$ minors of a $2 \times n$ matrix vanish then the matrix is a concatenation of $2 \times n_i$ matrices of rank 1, separated by zero columns.'' In fact, the number of associated primes of these ideals (for $2 \times n$ matrices) is exactly the $n$th Fibonacci number! If certain polynomial equations characterize conditions of interest, then the associated primes of the ideal they generate, obtained through primary decomposition, can help us understand the pieces.

\begin{figure}[!htb]
    \centering
    \includegraphics[width=0.35\textwidth]{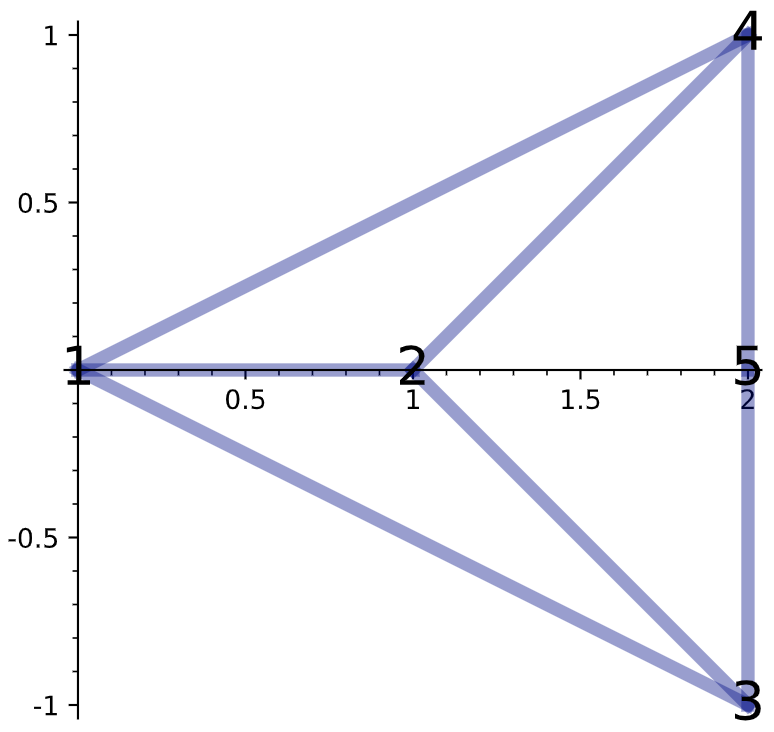}
    \caption{Slingshot example}
    \label{figure:slingshot-configuration}
\end{figure}

\begin{figure}[!htb]
    \centering
    \includegraphics[width=0.6\textwidth]{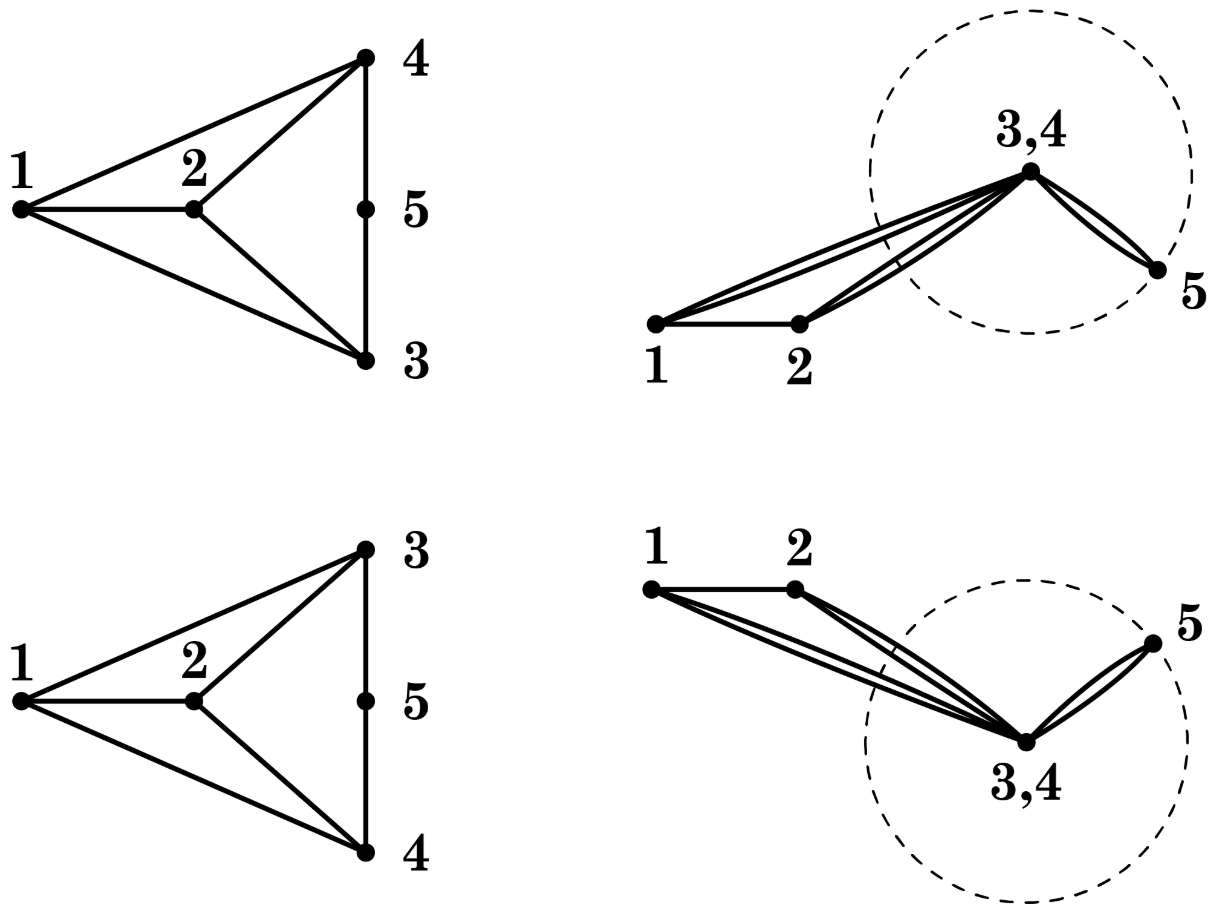}
    \caption{Some associated primes}
    \label{figure:four-slingshots}
\end{figure}

We now demonstrate this technique on a smaller bar framework which we call the \textit{slingshot}. Figure \ref{figure:slingshot-configuration} displays the graph with $V = [5], E =\{12,13,14,23,24,35,45\}$ in a specific configuration, but instead we will consider \textit{all possible configurations}.
We can calculate the rigidity matrix $dg$ after a change of coordinates to a \textit{moving frame}, as discussed in Section \ref{section:numerical-algebraic-geometry}, obtaining:
\begin{equation*}\footnotesize
    \left(\begin{array}{rrrrrrrrrr}
-x_{21} & 0 & x_{21} & 0 & 0 & 0 & 0 & 0 & 0 & 0 \\
-x_{31} & -x_{32} & 0 & 0 & x_{31} & x_{32} & 0 & 0 & 0 & 0 \\
-x_{41} & -x_{42} & 0 & 0 & 0 & 0 & x_{41} & x_{42} & 0 & 0 \\
0 & 0 & x_{21} - x_{31} & -x_{32} & -x_{21} + x_{31} & x_{32} & 0 & 0 & 0 & 0 \\
0 & 0 & x_{21} - x_{41} & -x_{42} & 0 & 0 & -x_{21} + x_{41} & x_{42} & 0 & 0 \\
0 & 0 & 0 & 0 & x_{31} - x_{51} & x_{32} - x_{52} & 0 & 0 & -x_{31} + x_{51} & -x_{32} + x_{52} \\
0 & 0 & 0 & 0 & 0 & 0 & x_{41} - x_{51} & x_{42} - x_{52} & -x_{41} + x_{51} & -x_{42} + x_{52}
\end{array}\right)
\end{equation*}
The generic rank of this matrix is $7$, so we will consider the $120$ minors of size $7$. Only $95$ of these minors are nonzero polynomials. Here is one of them:
\begin{multline*}
    x_{21}^{2} x_{32}^{2} x_{41} x_{42}^{2} - x_{21}^{2} x_{31} x_{32} x_{42}^{3} - x_{21}^{2} x_{32}^{2} x_{42}^{2} x_{51} + x_{21}^{2} x_{32} x_{42}^{3} x_{51} - x_{21}^{2} x_{32}^{2} x_{41} x_{42} x_{52}\\ + 2 \, x_{21}^{2} x_{31} x_{32} x_{42}^{2} x_{52} - x_{21}^{2} x_{32} x_{41} x_{42}^{2} x_{52} + x_{21}^{2} x_{32}^{2} x_{42} x_{51} x_{52}\\ - x_{21}^{2} x_{32} x_{42}^{2} x_{51} x_{52} - x_{21}^{2} x_{31} x_{32} x_{42} x_{52}^{2} + x_{21}^{2} x_{32} x_{41} x_{42} x_{52}^{2}
\end{multline*}
We also include the $7$ member constraints (\ref{equation:member-constraints}), which are:
\begin{equation*}
    \begin{array}{c}
        x_{21}^{2} - 1\\
x_{31}^{2} + x_{32}^{2} - 5\\
x_{41}^{2} + x_{42}^{2} - 5\\
{\left(x_{21} - x_{31}\right)}^{2} + x_{32}^{2} - 2\\
{\left(x_{21} - x_{41}\right)}^{2} + x_{42}^{2} - 2\\
{\left(x_{31} - x_{51}\right)}^{2} + {\left(x_{32} - x_{52}\right)}^{2} - 1\\
{\left(x_{41} - x_{51}\right)}^{2} + {\left(x_{42} - x_{52}\right)}^{2} - 1
    \end{array}.
\end{equation*}
This gives a total of $102$ equations in $nd - \binom{d+1}{2} = 9$ variables. However, code available at \cite{heaton-tensegrity-trusses-github} collects these $102$ equations in a list named \texttt{eqnz} and then calculates
\begin{verbatim}
xvarz = [var('x%s%s'%(i,k)) for i in range(1,6+1) for k in range(1,2+1) if i > k]
R = PolynomialRing(QQ,xvarz)
eqnz = [R(eqn) for eqn in eqnz]
I = R.ideal(eqnz)
AP = I.associated_primes()
\end{verbatim}
The output of this calculation is $8$ associated prime ideals, each described by a list of its generating polynomials.
\begin{verbatim}
[x52, x51 + 2, x42 - 1, x41 + 2, x32 + 1, x31 + 2, x21 + 1]
[x52, x51 + 2, x42 + 1, x41 + 2, x32 - 1, x31 + 2, x21 + 1]
[x42 + 1, x41 + 2, x32 + 1, x31 + 2, x21 + 1, x51^2 + x52^2 + 4*x51 + 2*x52 + 4]
[x42 - 1, x41 + 2, x32 - 1, x31 + 2, x21 + 1, x51^2 + x52^2 + 4*x51 - 2*x52 + 4]
[x52, x51 - 2, x42 - 1, x41 - 2, x32 + 1, x31 - 2, x21 - 1]
[x52, x51 - 2, x42 + 1, x41 - 2, x32 - 1, x31 - 2, x21 - 1]
[x42 + 1, x41 - 2, x32 + 1, x31 - 2, x21 - 1, x51^2 + x52^2 - 4*x51 + 2*x52 + 4]
[x42 - 1, x41 - 2, x32 - 1, x31 - 2, x21 - 1, x51^2 + x52^2 - 4*x51 - 2*x52 + 4]
\end{verbatim}
While trying to understand $102$ equations by inspection is too difficult, we can understand the output above. First we notice that $x_{21}$ is either $1$ or $-1$ (which we actually could have predicted from the $102$ equations). We could also notice that $x_{32} \in \{1,-1\}$ as well, but then our illustration would not explicitly show reflection symmetry. In conclusion, we draw illustrations of exactly $4$ of the eight prime ideals in Figure \ref{figure:four-slingshots}.

\textbf{Acknowledgements:} The author would like to thank Andrew Frohmader, Istvan Lauko, and Gabriella Pinter for discussions on applied math throughout the Spring of 2019, Timothy Duff, Sascha Timme, and Paul Breiding for education on numerical algebraic geometry, Robert Connelly, Miranda Holmes-Cerfon, and Louis Theran for discussions on rigidity theory, and Gabriella Pinter and Bernd Sturmfels for encouragement to write this article. Further, the author would like to thank Myfawny Evans, Frank Lutz, and Bernd Sturmfels for the opportunity to study rigidity theory as part of the Math+ Einstein Thematic Semester on Geometric and Topological Structure of Materials.

\bibliographystyle{plain}
\bibliography{references}

\begin{footnotesize}
Alexander Heaton (Max Planck Institute for Mathematics in the Sciences and Technische Universit\"at Berlin) can be reached at alexheaton2@gmail.com, or heaton@mis.mpg.de.
\end{footnotesize}

\end{document}